\documentclass[12pt]{amsart}

\usepackage{amsmath, amsthm, amscd, amsfonts, amssymb, graphicx, color}
\usepackage{times}
\usepackage{bm}
\newtheorem{theorem}{Theorem}

\def\ca{\mathring}

\newtheorem{lemma}{Lemma}
\DeclareSymbolFont{AMSb}{U}{msb}{m}{n}
\DeclareMathSymbol{\N}{\mathbin}{AMSb}{"4E}
\DeclareMathSymbol{\Z}{\mathbin}{AMSb}{"5A}
\DeclareMathSymbol{\R}{\mathbin}{AMSb}{"52}
\DeclareMathSymbol{\Q}{\mathbin}{AMSb}{"51}
\DeclareMathSymbol{\I}{\mathbin}{AMSb}{"49}
\DeclareMathSymbol{\C}{\mathbin}{AMSb}{"43}
\def\ca{\mathring}
\def\endex{{\hfill$\square$\medskip}}

\begin{document}

\title{Stability and uniqueness of $p$-values for likelihood-based inference}

\author{THOMAS J. DICICCIO}
\address{Department of Social Statistics,
Cornell University, Ithaca, New York 14853, U.S.A.} \email{tjd9@cornell.edu}

\author{TODD A. KUFFNER}
\address{Department of Mathematics, Washington University in St. Louis, St. Louis, Missouri 63130, U.S.A.} \email{kuffner@math.wustl.edu}

\author{G. ALASTAIR YOUNG}
\address{Department of Mathematics, Imperial College London, London SW7 2AZ,
U.K.} \email{alastair.young@imperial.ac.uk}

\author{RUSSELL ZARETZKI}
\address{Department of Statistics, University of Tennessee, Knoxville, Knoxville, Tennessee 37996, U.S.A.}
\email{rzaretzki@utk.edu}

\begin{abstract}
Likelihood-based methods of statistical inference provide a useful general methodology that is appealing, as a straightforward asymptotic theory can be applied for their implementation. It is important to assess the relationships between different likelihood-based inferential procedures in terms of accuracy and adherence to key principles of statistical inference, in particular those relating to conditioning on relevant ancillary statistics.  An analysis is given of the stability properties of a general class of likelihood-based statistics, including those derived from forms of adjusted profile likelihood, and comparisons are made between inferences derived from different statistics. In particular, we derive a set of sufficient conditions for agreement to $O_{p}(n^{-1})$, in terms of the sample size $n$, of inferences, specifically $p$-values, derived from different asymptotically standard normal pivots. Our analysis includes inference problems concerning a scalar or vector interest parameter, in the presence of a nuisance parameter. 
\end{abstract}


\maketitle

\section{Introduction}

A highly useful statistical methodology for inference on a scalar or vector interest parameter in the presence of a nuisance parameter is furnished by procedures based on the likelihood function, including tests and confidence sets based on the likelihood ratio statistic. Though no explicit optimality criteria are invoked, a quite general asymptotic theory allows straightforward implementation of such methodology in a wide range of settings. However,
accuracy and what may be termed inferential correctness are (Young (2009)) key desiderata of any parametric inference. When constructing, say, a confidence set for a parameter of interest in the presence of nuisance parameters, we desire high levels of coverage accuracy from the confidence set. Further, it is important that procedures are inferentially correct, meaning that they respect key principles of inference, in particular those relating to appropriate conditioning on ancillary information when this is relevant. The crucial issue here is the stability of the statistic used for inference, the extent to which the unconditional distribution of the statistic agrees with the conditional distribution of the statistic, relevant for achieving inferential correctness. Henceforth, when speaking of the stability of a pivot, we mean whether or not its marginal distribution inherently respects ancillary information. Specifically, a statistic which is stable to second-order is one whose conditional distribution given the observed value of an ancillary statistic agrees to second-order, $O(n^{-1})$, in the sample size $n$ with its marginal distribution. Our objective in this paper is to both analyse and elucidate properties of likelihood-based methods of statistical inference against these desiderata, and to provide new results that shed light on what is achieved by alternative approaches to implementation of likelihood-based methods of inference. We make two novel contributions.

We provide a general assessment of the stability properties of likelihood-based statistics commonly used for parametric inference. Our analysis considers first the case of the signed root likelihood ratio statistic for inference on a scalar interest parameter, in the presence of a nuisance parameter. In doing so, we establish a generalization to the practically realistic context involving nuisance parameters of results described by McCullagh (1984) and Severini (1990). We then discuss this issue for asymptotically standard normal pivots more generally, in particular those constructed from adjusted forms of profile likelihood, before considering inference for vector interest parameters. The results presented here allow comparisons to be drawn between the inferential properties of parametric bootstrap procedures and techniques of higher-order inference based on asymptotic, analytic approximation.

We also provide an explicit comparison of inferences, specifically $p$-values, obtained from different asymptotically standard normal pivots, including those constructed from adjusted forms of profile likelihood, establishing certain higher-order equivalences and differences. We derive a set of sufficient conditions ensuring agreement of $p$-values derived from different asymptotically standard normal pivots, to order $O_{p}(n^{-1})$.

\section{Background}

Suppose that $Y=(Y_1,\ldots,Y_n)$ is a continuous random vector and that the distribution of $Y$ depends on an unknown $d$-dimensional parameter $\theta$, partitioned as $\theta=(\psi,\phi)$, where initially we suppose $\psi=\theta_{1}$ is a scalar interest parameter and $\phi$ is a nuisance parameter of dimension $d-1$. We later consider the case of a vector interest parameter $\psi$.

Let $L(\theta)$ be the loglikelihood function for $\theta$ based on $Y$ and let $\hat\theta=(\hat\psi,\hat\phi)$ be the global maximum likelihood estimator of $\theta$.
Further, let $\tilde\theta=\tilde\theta(\psi)=(\psi,\tilde\phi)=\{\psi,\tilde\phi(\psi)\}$ be the constrained maximum likelihood estimator of $\theta$ for given $\psi$. Then the profile loglikelihood function for $\psi$ is $M(\psi)=L\{\tilde\theta(\psi)\}$ and the likelihood ratio statistic for $\psi$ is $W(\psi)=2\{M(\hat\psi)-M(\psi)\}$, where $M(\hat\psi)=L(\hat\theta)$, since $\tilde\theta(\hat\psi)=\hat\theta$. The signed root likelihood ratio statistic is $R(\psi)={\rm sgn}(\hat\psi-\psi)\{W(\psi)\}^{1/2}$. Testing $H_0:\psi=\psi_0$ against $H_a:\psi>\psi_0$ or $H_a:\psi<\psi_0$ can be based on the test statistic $R(\psi_0)$. Asymptotically, as the sample size $n$ increases, the sampling distribution of $R(\psi)$ tends to the standard normal distribution. Heading the list of desiderata for refinement of the inference procedures furnished by such first-order asymptotic theory is the achievement of higher-order accuracy in distributional approximation, while respecting the need for inferential correctness.

Two main routes (Young (2009)) to higher-order accuracy emerge from contemporary statistical theory. The most developed route is that which utilises analytic procedures, based on `small-sample asymptotics', such as saddlepoint approximation and related methods, to refine first-order distribution theory. The second route involves simulation or bootstrap methods, which aim to obtain refined distributional approximations directly, without analytic approximation: see, for instance, DiCiccio, Martin and Stern (2001), Lee and Young (2005), DiCiccio and Young (2008).

A detailed account of analytic methods for distributional approximation which yield higher-order accuracy is given by Barndorff-Nielsen and Cox (1994). Two particular highlights of an intricate theory are especially important:   Bartlett correction of the likelihood ratio statistic $W(\psi)$, which we discuss in Section 8, and the construction of analytically modified forms of the signed root likelihood ratio statistic $R(\psi)$, designed to offer higher-order accuracy. These procedures also provide inferential correctness, specifically conditional validity, to high (asymptotic) order, in the two key settings where conditional inference is crucial,  namely multi-parameter exponential family and ancillary statistic contexts. Particularly central to the analytic approach to higher-order accurate inference on a scalar interest parameter is Barndorff-Nielsen's $R^*$ statistic (Barndorff-Nielsen (1986)). In both the multi-parameter exponential family and ancillary statistic contexts, the $R^*$ statistic is conditionally, and hence unconditionally, distributed as standard normal, to error of third-order $O(n^{-3/2})$ in the sample size. So, analytic standard normal approximation of the sampling distribution of the $R^*$ statistic yields third-order accuracy under repeated sampling, while respecting the requirements of conditioning to that same order.

Lawley (1956) showed that $E_{\theta}\{R(\psi)\}=n^{-1/2}m(\theta)+O(n^{-3/2})$ and ${\rm var}_{\theta}\{R(\theta)\}=1+n^{-1}v(\theta)+O(n^{-2})$, where $m(\theta)$ and $v(\theta)$ are both of order $O(1)$, while the third and higher-order cumulants are of order $O(n^{-3/2})$ or smaller; see also Bickel and Ghosh (1990).
Therefore, $\{R(\psi)-n^{-1/2}m(\theta)\}/\{1+n^{-1}v(\theta)\}^{1/2}$ has the standard normal distribution to error of order $O(n^{-3/2})$.
DiCiccio and Stern (1994a) showed that $\{R(\psi)-n^{-1/2}m(\tilde\theta)\}/\{1+n^{-1}v(\tilde\theta)\}^{1/2}$ also has the standard normal distribution to error of order $O(n^{-3/2})$.
This DiCiccio and Stern (1994a) result asserts that $[R(\psi)-E_{\tilde\theta}\{R(\psi)\}]/{[{\rm var}_{\tilde\theta}\{R(\psi)\}]^{1/2}}$ is also distributed as standard normal to error of order $O(n^{-3/2})$. In turn, this distributional result immediately suggests the parametric bootstrap approaches to third-order accurate inference discussed by DiCiccio et al. (2001) and Lee and Young (2005). For testing $H_0:\psi=\psi_0$ against one-sided alternatives, $p-$values distributed, under repeated sampling, as uniform to error of order $O(n^{-3/2})$, and hence yielding error rate $O(n^{-3/2})$, can be obtained by bootstrapping $R(\psi_0)$ at the parameter value $\theta=(\psi_0,\tilde\phi_0)$,  where $\tilde\phi_0=\tilde\phi(\psi_0)$. DiCiccio and Young (2008) show that this parametric bootstrap procedure respects the requirements of conditioning in multi-parameter exponential family settings to third-order.

From a repeated sampling perspective, such third-order accurate inference can be similarly obtained (Lee and Young, 2005) by bootstrap approximation to the sampling distribution of other asymptotically standard normal pivots, in particular, pivots constructed as standardized versions of the difference $\hat \psi-\psi_0$ or the score function $\partial M(\psi)/\partial\psi |_{\psi=\psi_0}$, that avoid calculation of both the global and constrained maximum likelihood estimators, and may therefore may be more appealing for use in a computationally-intensive bootstrap inference. A fundamental question that arises concerns the inferential implications of choice of a particular statistic: when do inferences based on different choices of statistic agree to high-order? It is also necessary to ask whether such inference respects the requirements of conditioning on relevant ancillary statistics, in models which admit the existence of such. Since a bootstrap calculation involves unconditional sampling at parameter value $\theta=(\psi_0,\tilde\phi_0)$, the key question is the extent to which the conditional and unconditional distributions of the statistic being used for the inference differ.

In this paper we provide an analysis directed at these questions, providing new results on the stability properties of likelihood-based statistics and agreement of $p$-values derived from different asymptotically normal pivots. The implications of the analysis for bootstrap methodology and detailed comparisons of the latter with analytic procedures of inference will be described elsewhere.

We consider first the stability properties of the signed root statistic $R(\psi)$; in doing so, we establish a generalization to the nuisance parameter context of a result of McCullagh (1984): see also Severini (2000, Section 6.4.4).
We then discuss the stability issue in problems involving nuisance parameters for asymptotically standard normal pivots more generally, before examining conditions which ensure that $p$-values derived from two different pivots agree to second-order. Extension of the conclusions to test statistics based on general adjusted forms of profile likelihood are described, before presenting results concerning inference for vector interest parameters.

Our analysis is concerned exclusively with inferential comparisons `under the null' so, for instance we examine the unconditional and conditional distributions of the signed root statistic $R(\psi)$ under the model in question when the true parameter value is $\theta=(\psi, \phi)$. Similarly, the analysis concerns comparison of different $p$-values under assumed correctness of the null hypothesis being tested.

\section{Notation}

In the calculations that follow, arrays and summation are denoted by using the standard conventions, for which the indices $r,s,t,\ldots$ are assumed to range over $1,\ldots,d$. Summation over the range is implied for any index appearing in an expression both as a subscript and as a superscript.
Differentiation is indicated by subscripts, so $L_r(\theta)=\partial L(\theta)/\partial\theta^r$, $L_{rs}(\theta)=\partial^2 L(\theta)/\partial\theta^r\partial\theta^s$, etc.
Then $E\{L_r(\theta)\}=0$; let $\lambda_{rs}=E\{L_{rs}(\theta)\}$, $\lambda_{rst}=E\{L_{rst}(\theta)\}$, etc., and
put $l_r=L_r(\theta)$, $l_{rs}=L_{rs}(\theta)-\lambda_{rs}$, $l_{rst}=L_{rst}(\theta)-\lambda_{rst}$, etc.
The constants $\lambda_{rs}$, $\lambda_{rst}, \ldots$, are assumed to be of order $O(n)$.
The variables $l_r$, $l_{rs}$, $l_{rst}$, etc., each of which have expectation 0, are assumed to be of order $O_p(n^{1/2})$.
The joint cumulants of $l_r$, $l_{rs}$, etc.\ are assumed to be of order $O(n)$. These assumptions are usually satisfied in situations involving independent
observations. The observed information matrix is $J(\theta)=[-L_{rs}(\theta)]$, while the expected (Fisher) information matrix is $I(\theta)=[-\lambda_{rs}(\theta)]$.
It is useful to extend the $\lambda$-notation: let $\lambda_{r,s}=E(L_rL_s)=E(l_rl_s)$, $\lambda_{rs,t}=E(L_{rs}L_t)=E(l_{rs}l_t)$, etc.
The Bartlett identities involving the $\lambda$'s can be derived by repeated differentiation of the identity $\int\exp\{L(\theta)\}dy=1$; in particular,
$$
\lambda_{rs}+\lambda_{r,s}=0, \quad
\lambda_{rst}+\lambda_{rs,t}+\lambda_{rt,s}+\lambda_{st,r}+\lambda_{r,s,t}=0.
$$
Differentiation of the definition $\lambda_{rs}=\int L_{rs}(\theta)\exp\{L(\theta)\}dy$ yields
$\lambda_{rs/t}=\lambda_{rst}+\lambda_{rs,t},$
where $\lambda_{rs/t}=\partial\lambda_{rs}/\partial\theta^t$.
Further, let $(\lambda^{rs})$ be the $d \times d$ matrix inverse of $(\lambda_{rs})$, and let $\eta=-1/\lambda^{11}$, $\tau^{rs}=\eta\lambda^{1r}\lambda^{1s}$,  and $\nu^{rs}=\lambda^{rs}+\tau^{rs}$.
Thus, $\lambda^{rs}$, $\tau^{rs}$, and $\nu^{rs}$ are of order $O(n^{-1})$, while $\eta$ is of order $O(n)$. For clarity, we point out that a superscript or subscript of `$1$' refers to the scalar interest parameter $\psi$, where $\psi$ is the first component of $\theta$.

Suppose that $A$ is an ancillary, i.e., distribution constant, statistic such that $(\hat\theta,A)$ is sufficient.
To distinguish conditional calculations from unconditional ones, the accent symbol $\mathring{~~~}$ is used to denote quantities derived from the conditional distribution of $Y$ given $A$.
Since the conditional loglikelihood $\mathring{L}(\theta)$ differs from the unconditional loglikelihood $L(\theta)$ by a quantity that depends on $A$ but not on $\theta$, it follows that $\mathring{W}(\psi)=W(\psi)$ and that $\mathring{L}_r=L_r$, $\mathring{L}_{rs}=L_{rs}$, etc.
Let $\mathring{\lambda}_{rs}=\mathring {E}\{L_{rs}(\theta)\}$, $\mathring{\lambda}_{rst}=\mathring{E}\{L_{rst}(\theta)\}$, etc.,
and put $\mathring{l}_r=l_r(\theta)$, $\mathring{l}_{rs}=L_{rs}(\theta)-\mathring{\lambda}_{rs}$, $\mathring{l}_{rst}=L_{rst}(\theta)-\mathring{\lambda}_{rst}$, etc.
The quantities $\mathring{\lambda}_{rs}$, $\mathring{\lambda}_{rst}$, etc.\ are random variables depending on $A$, assumed to be of order $O_p(n)$.
The variables $\mathring{l}_r$, $\mathring{l}_{rs}$, $\mathring{l}_{rst}$, etc.\ have conditional expectation 0, so they also have unconditional expectation 0, and they are assumed to be of order $O_p(n^{1/2})$.
Further, the joint conditional cumulants of $\mathring{l}_r, \mathring{l}_{rs}$, etc.\ depend on $A$, and they are assumed to be of order $O_p(n)$.
It is useful to extend the $\mathring{\lambda}$-notation by letting $\mathring{\lambda}_{r,s}=\mathring{E}(L_rL_s)=\mathring{E}(l_rl_s)$, $\mathring{\lambda}_{rs,t}=\mathring {E}(L_{rs}L_t)=\mathring {E}(l_{rs}l_t)$, etc.
Also, let $(\mathring\lambda^{rs})$ be the $d \times d$ matrix inverse of $(\mathring\lambda_{rs})$, and let $\mathring\eta=-1/\mathring\lambda^{11}$, $\mathring\tau^{rs}=\mathring\eta\mathring\lambda^{1r}\mathring\lambda^{1s}$,  and $\mathring\nu^{rs}=\mathring\lambda^{rs}+\mathring\tau^{rs}$,
so that $\mathring\lambda^{rs}$, $\mathring\tau^{rs}$, and $\mathring\nu^{rs}$ are of order $O_p(n^{-1})$, while $\mathring\eta$ is of order $O_p(n)$.

Following Barndorff-Nielsen and Cox (1994, Section 7.2), construction of an ancillary statistic $A$ such that $(\hat \theta, A)$ is sufficient is, except in rather special cases, only possible for transformation models and, in a degenerate sense, for full exponential family models, where $\hat \theta$ itself is sufficient. It is therefore in general necessary to consider conditioning on statistics $A$ which are approximately ancillary in a suitable sense. Results presented here continue to hold under the assumption that $A$ is locally ancillary (Cox (1980)). Let $\theta_0$ be an arbitrary but specified parameter value, and let $A \equiv A(Y, \theta_0)$ be a candidate ancillary statistic. If the density of $A$ under parameter value $\theta_0 +n^{-1/2}\delta$ satisfies
\[
f_A(a; \theta_0+n^{-1/2}\delta)=f_A(a;\theta_0)\{1+O(n^{-q/2})\},
\]
then (Cox (1980), McCullagh (1987, Section 8.3)) $A$ is said to be $q$-th order local ancillary in the vicinity of $\theta_0$. Note that this definition applies only to parameter values in an $O(n^{-1/2})$ neighbourhood of $\theta_0$: if $\theta_0$ is the true parameter value, as $n$ increases the likelihood function becomes negligible outside this neighbourhood. The loglikelihood function based on $A$ satisfies $L_A(\theta_0+n^{-1/2}\delta)=L_A(\theta_0)+O(n^{-q/2})$. As is the case in the no nuisance parameter context considered by Severini (1990) and McCullagh (1987, Section 8.4), results in Section 4 relating to stability of asymptotically standard normal pivots continue to hold for any second-order local ancillary $A$, as do results in Section 8 concerning stability of an adjusted profile likelihood ratio statistic. Essentially, the assumption of a second-order local ancillary is sufficient to ensure the relationships detailed below between conditional and unconditional cumulants.

The technique of proof used here to compare the conditional and unconditional distributions of asymptotically standard normal pivots to second order is a generalization of that described by Severini (2000, Chapter 6) in the case of a scalar interest parameter without nuisance parameters. For this technique, it is essential to compare the $\mathring\lambda$-quantities with their $\lambda$-counterparts.

We first investigate the difference between $\mathring{\lambda}_{rs}$ and $\lambda_{rs}$; note that $\lambda_{rs}=E(L_{rs})=E\{\mathring{E}(L_{rs})\}=E(\mathring{\lambda}_{rs})$.
Furthermore, ${\rm var}(\mathring{\lambda}_{rs})={\rm var}\{\mathring{E}(L_{rs})\} = {\rm var}(L_{rs})-E\{\mathring{\rm var}(L_{rs})\} = O(n)-E\{O_p(n)\} =O(n)$,
and consequently, $\mathring{\lambda}_{rs}=\lambda_{rs}+O_p(n^{1/2})$.
An identical argument shows that $\mathring{\lambda}_{rst}=\lambda_{rst}+O_p(n^{1/2})$, etc.

Assume that differentiation of the identity $\ca\lambda_{rs}=\lambda_{rs}+O_p(n^{1/2})$ yields $\ca\lambda_{rs/t}=\lambda_{rs/t}+O_p(n^{1/2})$, where $\mathring\lambda_{rs/t}=\partial\mathring\lambda_{rs}/\partial\theta^t$ and, as before, $\lambda_{rs/t}=\partial\lambda_{rs}/\partial\theta^t$. We note that, as a rule, differentiation of an asymptotic
relation will preserve the asymptotic order, but that care is necessary; see Barndorff-Nielsen and Cox (1994, Exercise 5.4) and Pace and Salvan (1994).  The asymptotic order of the difference between $\ca\lambda_{rs/t}$ and $\lambda_{rs/t}$ indicated here, therefore, actually constitutes an additional assumption of our calculations.
The preceding results imply $\mathring{\lambda}_{rs,t}=\lambda_{rs,t}+O_p(n^{1/2})$, since the Bartlett identities $\mathring{\lambda}_{rs/t}=\mathring{\lambda}_{rst}+\mathring{\lambda}_{rs,t}$ and ${\lambda}_{rs/t}={\lambda}_{rst}+{\lambda}_{rs,t}$ yield $\mathring{\lambda}_{rs,t}=\mathring{\lambda}_{rs/t}-\mathring{\lambda}_{rst}
=\lambda_{rs/t}-\lambda_{rst}+O_p(n^{1/2})=\lambda_{rs,t}+O_p(n^{1/2})$.
Define $\mathring{\Delta}_{rs}=\mathring{\lambda}_{rs}-\lambda_{rs}$, so that $\mathring{\Delta}_{rs}$ is a function of $\theta$ and $A$, having order $O_p(n^{1/2})$.
Then ${l}_{rs}=L_{rs}-{\lambda}_{rs}=(L_{rs}-\mathring\lambda_{rs})+(\mathring{\lambda}_{rs}-\lambda_{rs})=\mathring l_{rs}+\mathring{\Delta}_{rs}$.

\section{Stability result for $R(\psi)$ and other pivots}

We now consider the stability of $R(\psi)$ and other asymptotically standard normal pivots.

\subsection{$R(\psi)$ is a stable pivot to second order}
\begin{theorem}
The conditional and unconditional distributions of $R(\psi)$ agree to error of order $O(n^{-1})$, given the ancillary statistic $A$.
\end{theorem}
\emph{Proof.}
To error of order $O(n^{-1})$, the variance of $R(\psi)$ is $1$ and the third- and higher-order cumulants are $0$; the mean is of order $O(n^{-1/2})$.
The conditional distribution given $A$ has the same cumulant structure as the unconditional distribution. Thus, to show that the conditional and unconditional distributions agree to second-order, it suffices to show that $\mathring E\{R(\psi)\}=E\{R(\psi)\}+O_p(n^{-1})$.

Standard calculations, such as those given by Lawley (1956) and detailed in the Appendix of DiCiccio and Stern (1994b), show that $W(\psi)$ has the expansion
\begin{align*}
W(\psi)&=\tau^{rs}l_rl_s-2\lambda^{rt}\tau^{su}l_{rs}l_tl_u-\tau^{rt}\tau^{su}l_{rs}l_tl_u+\lambda^{ru}\nu^{sv}\tau^{tw}\lambda_{rst}l_ul_vl_w\\
&\hspace{20pt}+{\textstyle{1\over 3}}\tau^{ru}\tau^{sv}\tau^{tw}\lambda_{rst}l_ul_vl_w+O_p(n^{-1}).
\end{align*}

DiCiccio and Stern (1994b) showed that $R(\psi)$ may be decomposed as $R(\psi)=\eta^{1/2}\{R_1+R_2+O_p(n^{-3/2})\}$, where $R_1=-\lambda^{1r}l_r$ and
\[
R_2=\lambda^{1r}\lambda^{st}l_{rs}l_t+{\textstyle{1\over
2}}\lambda^{1r}\tau^{st}l_{rs}l_t -{\textstyle{1\over
2}}\lambda^{1r}\lambda^{su}\nu^{tv}\lambda_{rst}l_ul_v
-{\textstyle{1\over 6}}\lambda^{1r}\tau^{su}\tau^{tv}\lambda_{rst}l_ul_v.
\]
Here $R_1$ is of order $O_p(n^{-1/2})$ and $R_2$ is of order $O_p(n^{-1})$.
Since $E(R_1)=0$, it follows that
\[
E\{R(\psi)\}=\eta^{1/2}\{\lambda^{1r}\lambda^{st}\lambda_{rs,t}
+{\textstyle{1\over 2}}\lambda^{1r}\tau^{st}\lambda_{rs,t}
+{\textstyle{1\over 2}}\lambda^{1r}\lambda^{st}\lambda_{rst}
+{\textstyle{1\over 3}}\lambda^{1r}\tau^{st}\lambda_{rst}\}+O(n^{-1}).
\]
Note also that $R_1=-\lambda^{1r}l_r=-\lambda^{1r}\mathring l_r$ and
\begin{align*}
R_2&=\lambda^{1r}\lambda^{st}l_{rs}l_t
+{\textstyle{1\over 2}}\lambda^{1r}\tau^{st}l_{rs}l_t
-{\textstyle{1\over 2}}\lambda^{1r}\lambda^{su}\nu^{tv}\lambda_{rst}l_ul_v
-{\textstyle{1\over 6}}\lambda^{1r}\tau^{su}\tau^{tv}\lambda_{rst}l_ul_v\\
&=\lambda^{1r}\lambda^{st}\mathring l_{rs} \mathring l_t+\lambda^{1r}\lambda^{st}\mathring{\Delta}_{rs} \mathring l_t
+{\textstyle{1\over 2}}\lambda^{1r}\tau^{st} \mathring l_{rs} \mathring l_t
+{\textstyle{1\over 2}}\lambda^{1r}\tau^{st} \mathring{\Delta}_{rs} \mathring l_t\\
&\hspace{20pt}-{\textstyle{1\over 2}}\lambda^{1r}\lambda^{su}\nu^{tv}\lambda_{rst}\mathring l_u \mathring l_v
-{\textstyle{1\over 6}}\lambda^{1r}\tau^{su}\tau^{tv}\lambda_{rst}\mathring l_u \mathring l_v.
\end{align*}
Thus, since $\mathring E(R_1)=0$,
\begin{align*}
\mathring E\{R(\psi)\}&=\eta^{1/2}\{\lambda^{1r}\lambda^{st}\mathring \lambda_{rs,t}
+{\textstyle{1\over 2}}\lambda^{1r}\tau^{st} \mathring \lambda_{rs,t}
+{\textstyle{1\over 2}}\lambda^{1r}\lambda^{su}\nu^{tv}\lambda_{rst}\mathring \lambda_{uv}
+{\textstyle{1\over 6}}\lambda^{1r}\tau^{su}\tau^{tv}\lambda_{rst}\mathring \lambda_{uv}+O_p(n^{-3/2})\}\\
&=\eta^{1/2}\{\lambda^{1r}\lambda^{st} \lambda_{rs,t}
+{\textstyle{1\over 2}}\lambda^{1r}\tau^{st}  \lambda_{rs,t}
+{\textstyle{1\over 2}}\lambda^{1r}\lambda^{su}\nu^{tv}\lambda_{rst} \lambda_{uv}
+{\textstyle{1\over 6}}\lambda^{1r}\tau^{su}\tau^{tv}\lambda_{rst} \lambda_{uv}+O_p(n^{-3/2})\}\\
&=\eta^{1/2}\{\lambda^{1r}\lambda^{st}\lambda_{rs,t}
+{\textstyle{1\over 2}}\lambda^{1r}\tau^{st}\lambda_{rs,t}
+{\textstyle{1\over 2}}\lambda^{1r}\lambda^{st}\lambda_{rst}
+{\textstyle{1\over 3}}\lambda^{1r}\tau^{st}\lambda_{rst}+O_p(n^{-3/2})\}\\
&=E\{R(\psi)\}+ O_p(n^{-1}).
\end{align*}
It follows that the conditional distribution of $R(\psi)$ differs from its marginal distribution by error of order $O(n^{-1})$, given $A$.
\endex

McCullagh (1984) generalized the notion of the signed root statistic to the case of a vector interest parameter and established this stability result in the case of no nuisance parameters; Severini (1990) gave a further demonstration for the case of a scalar interest parameter with no nuisance parameters.
Therefore, the result shown here extends the work of McCullagh and Severini to situations where nuisance parameters are present.

This second-order stability of $R(\psi)$ for the nuisance parameter context has been discussed, but not demonstrated formally as we have here, by Pierce and Bellio (2006).
The methodological consequence of the result is immediate. Any approximation to the unconditional distribution of $R(\psi)$ having error of order $O(n^{-1})$ also approximates the conditional distribution of $R(\psi)$ to the same order of error. Such an approximation may (DiCiccio et al. (2001)) be derived, for instance, from the bootstrap distribution of $R(\psi)$. If that approximation is then used, say, to construct confidence limits for $\psi$, then those limits have coverage error of order $O(n^{-1})$, conditionally as well as unconditionally.

\subsection{Stability of other asymptotically standard normal pivots}

We now consider general asymptotically standard normal pivots of the form $T(\psi)=\eta^{1/2}\{T_1+T_2+O_p(n^{-3/2})\}$, where $T_1=-\lambda^{1r}l_r$ and $T_2$ is of the form $T_2=\xi^{rst}l_{rs}l_t-\xi^{rs}l_rl_s$, with $\xi^{rst}$ and $\xi^{rs}$ assumed to be of order $O(n^{-2})$, so that $T_1$ is of order $O_p(n^{-1/2})$ and $T_2$ is of order $O_p(n^{-1})$.
 We demonstrate that commonly used pivots may all be expressed in this form; for example, for $R(\psi)$, the preceding expansions show that $\xi^{rst}=\lambda^{1r}\lambda^{st}+{\textstyle{1\over 2}}\lambda^{1r}\tau^{st}$ and $\xi^{rs}={\textstyle{1\over 2}}\lambda^{1t}\lambda^{ur}\nu^{vs}\lambda_{tuv}+{\textstyle{1\over 6}}\lambda^{1t}\tau^{ur}\tau^{vs}\lambda_{tuv}$.
Both conditionally and unconditionally, the fourth- and higher-order cumulants of such a pivot are immediately seen to be of order $O(n^{-1})$ or smaller. Consequently, if we are to show that the conditional and unconditional distributions of these pivots agree to error of order $O(n^{-1})$ given $A$, all we need to show is that the first three conditional cumulants agree with the unconditional ones to error of order $O_p(n^{-1})$.
We show that the first and third conditional cumulants agree with the unconditional ones to the required order of error without further restrictions on $\xi^{rs}$ and $\xi^{rst}$. We demonstrate that for the second conditional cumulant to agree to with the unconditional one a sufficient condition is that $\xi^{rs1}={\textstyle{1\over 2}}\lambda^{1r}\lambda^{1s}$.
It is easy to see that $R(\psi)$ satisfies this criterion for, in this case,
\begin{align*}
\xi^{rs1}=\lambda^{1r}\lambda^{s1}+{\textstyle{1\over 2}}(\lambda^{1r}\eta\lambda^{1s}\lambda^{11})
=\lambda^{1r}\lambda^{1s}+{\textstyle{1\over 2}}\{\lambda^{1r}(-1/\lambda^{11})\lambda^{1s}\lambda^{11}\}
=\lambda^{1r}\lambda^{1s}-{\textstyle{1\over 2}}\lambda^{1r}\lambda^{1s}={\textstyle{1\over 2}}\lambda^{1r}\lambda^{1s}.
\end{align*}
\begin{theorem}
The unconditional and conditional distributions of $T(\psi)$ agree to error of order $O(n^{-1})$ given the ancillary statistic $A$.
\end{theorem}
The result follows immediately from three lemmas concerning the stability of the first three cumulants of $T(\psi)$, beginning with the first cumulant, the mean.
\begin{lemma}
$\mathring E\{T(\psi)\}=E\{T(\psi)\}+O_p(n^{-1})$.
\end{lemma}
\emph{Proof.}
Recall that $T_1=-\lambda^{1r}l_r=-\lambda^{1r}\mathring l_r$ and that $T_2=\xi^{rst}l_{rs}l_t-\xi^{rs}l_rl_s=\xi^{rst}(\mathring l_{rs}+\mathring\Delta_{rs})\mathring l_t-\xi^{rs}\mathring l_r\mathring l_s$.
Then, $E\{T(\psi)\}=\eta^{1/2}\{\xi^{rst}\lambda_{rs,t}+\xi^{rs}\lambda_{rs}+O(n^{-3/2})\}$ and
\begin{align*}
\mathring E\{T(\psi)\}&=\eta^{1/2}\{\xi^{rst}\mathring\lambda_{rs,t}+\xi^{rs}\mathring\lambda_{rs}+O_p(n^{-3/2})\}\\
&=\eta^{1/2}\{\xi^{rst}\lambda_{rs,t}+\xi^{rs}\lambda_{rs}+O_p(n^{-3/2})\}.
\end{align*}
Therefore, the conditional first cumulant agrees with the unconditional one to error of order $O_p(n^{-1})$, as required.
\endex

\begin{lemma}
If $\xi^{rs1}={\textstyle{1\over 2}}\lambda^{1r}\lambda^{1s}$, then $\mathring{{\rm var}}\{T(\psi)\}={\rm var}\{T(\psi)\}+O_p(n^{-1})$.
\end{lemma}
\emph{Proof.}
See Appendix.
\endex

\begin{lemma}
$\mathring{{\rm skew}}\{T(\psi)\}={{\rm skew}}\{T(\psi)\}+O_p(n^{-1})$.
\end{lemma}
\emph{Proof.}
See Appendix.
\endex

A sufficient condition for $\mathring{{\rm var}}\{T(\psi)\}={\rm var}\{T(\psi)\}+O_p(n^{-1})$ is $\xi^{rs1}={\textstyle{1\over 2}}\lambda^{1r}\lambda^{1s}$; if this holds, we have ${\rm skew}\{T(\psi)\}=\eta^{3/2}(\lambda^{1r}\lambda^{1s}\lambda^{1t}\lambda_{rst}-6\xi^{11})+O(n^{-1})$.

\section{Comparison of $p$-values}

Our objective here is to utilize preceding calculations to examine conditions which ensure that $p$-values based on two different asymptotically normal pivots agree to second-order. Here we refer to the $p$-value calculated from the exact sampling distribution of the pivot, or any approximation to the exact $p$-value accurate to $O_p(n^{-1})$. Such accuracy of approximation is obtained, for instance, quite generally for an asymptotically normal pivot by bootstrapping (Lee and Young (2005)), but would not be obtained by the normal approximation.

Consider hypothesis testing for $\psi$ based on a test statistic expressible as  $T(\psi)=\eta^{1/2}(T_1+T_2)+O_p(n^{-1})$, where $T_1=-\lambda^{1r}l_r$ and $T_2$ is of the form $T_2=\xi^{rst}l_{rs}l_t-\xi^{rs}l_rl_s$, with $\xi^{rst}$ and $\xi^{rs}$ assumed to be of order $O(n^{-2})$.
We have shown that the first three cumulants of $T(\psi)$ are
\begin{align*}
\kappa_1&=E\{T(\psi)\}=\eta^{1/2}(\xi^{rst}\lambda_{rs,t}+\xi^{rs}\lambda_{rs})+O(n^{-1}),\\
\kappa_2&={\rm var}\{T(\psi)\}=1+O(n^{-1}),\\
\kappa_3&={\rm skew}\{T(\psi)\}=\eta^{3/2}(\lambda^{1r}\lambda^{1s}\lambda^{1t}
\lambda_{rst}+3\lambda^{1r}\lambda^{1s}\lambda^{1t}\lambda_{rs,t}-6\xi^{rs1}
\lambda^{1t}\lambda_{rs,t}-6\xi^{11})+O(n^{-1}),
\end{align*}
while the fourth- and higher-order cumulants are of order $O(n^{-1})$ or smaller.

Consider another test statistic $\breve T(\psi)=\eta^{1/2}(\breve T_1+\breve T_2)+O_p(n^{-1})$, where $\breve T_1=-\lambda^{1r}l_r=T_1$ and $\breve T_2$ is of the form $\breve T_2=\breve \xi^{rst}l_{rs}l_t-\breve \xi^{rs}l_rl_s$, with $\breve \xi^{rst}$ and $\breve \xi^{rs}$ assumed to be of order $O(n^{-2})$. Our goal is to establish conditions on the two pivots $T(\psi)$ and $\breve T(\psi)$ which ensure that $p$-values agree to second-order.

\begin{theorem}
If the conditions
\begin{equation}
\breve\xi^{rst}=\xi^{rst}+O(n^{-5/2}),\label{E:first condition}
\end{equation}
\begin{equation}
\breve\xi^{rs}+\breve\xi^{tu}\lambda_{tu}\tau^{rs}=\xi^{rs}+\xi^{tu}\lambda_{tu}\tau^{rs}+O(n^{-5/2}),\label{E:second condition}
\end{equation}
are satisfied, then the $p$-value derived from the pivot $T(\psi)$ agrees with that derived from the pivot $\breve T(\psi)$ to error of order $O_p(n^{-1})$.
\end{theorem}
\emph{Proof.}
The $p$-value for testing against alternatives greater than $\psi$ is the right-hand tail probability for $T(\psi)$.
The normalizing Cornish-Fisher expansion shows that the $p$-value is
\[
1-\Phi(\eta^{1/2}T_1+\eta^{1/2}T_2-{\textstyle{1\over 6}}\kappa_3\eta T_1^2-\kappa_1+{\textstyle{1\over 6}}\kappa_3)+O_p(n^{-1}),
\]
where $\Phi(\cdot)$ denotes the standard normal cumulative distribution function.

Let the first three cumulants of $\breve T(\psi)$ be denoted by $\breve \kappa_1$, $\breve \kappa_2$, $\breve \kappa_3$; the $p$-value based on $\breve T(\psi)$ is
\[
1-\Phi(\eta^{1/2}T_1+\eta^{1/2}\breve T_2-{\textstyle{1\over 6}}\breve \kappa_3\eta T_1^2-\breve \kappa_1+{\textstyle{1\over 6}}\breve \kappa_3)+O_p(n^{-1}).
\]
We now determine sufficient conditions on $\breve \xi^{rs}$ and $\breve \xi^{rst}$ to ensure that the $p$-value obtained from $\breve T(\psi)$ agrees with that obtained from $T(\psi)$ to error of order $O_p(n^{-1})$.
Agreement of the $p$-values to this order occurs when
\[
\eta^{1/2}\breve T_2-{\textstyle{1\over 6}}\breve \kappa_3\eta T_1^2-\breve \kappa_1+{\textstyle{1\over 6}}\breve \kappa_3
=\eta^{1/2}T_2-{\textstyle{1\over 6}}\kappa_3\eta T_1^2-\kappa_1+{\textstyle{1\over 6}}\kappa_3
\]
to error of order $O_p(n^{-1})$, that is when
\[
\{\eta^{1/2}(\breve T_2-T_2)-{\textstyle{1\over 6}}(\breve \kappa_3-\kappa_3)\eta T_1^2\}
-\{(\breve \kappa_1-\kappa_1)-{\textstyle{1\over 6}}(\breve \kappa_3-\kappa_3)\}=O_p(n^{-1}).
\]
The first term on the left-hand side of the preceding equation is random, as it involves terms of the form $l_{rs}l_t$ and $l_rl_t$, while the second term is a constant. Consequently, by separating the random and non-random components, we see that the preceding equation actually stipulates two conditions:
\begin{align*}
\eta^{1/2}(\breve T_2-T_2)-{\textstyle{1\over 6}}(\breve \kappa_3-\kappa_3)\eta T_1^2&=O_p(n^{-1}),\\
(\breve \kappa_1-\kappa_1)-{\textstyle{1\over 6}}(\breve \kappa_3-\kappa_3)&=O(n^{-1}).
\end{align*}
The second of these equations gives $(\breve \kappa_1-\kappa_1)={\textstyle{1\over 6}}(\breve \kappa_3-\kappa_3)+O(n^{-1})$,
so we can write the equations as
\begin{align}
\eta^{1/2}(\breve T_2-T_2)-(\breve \kappa_1-\kappa_1)\eta T_1^2&=O_p(n^{-1}),\label{E:first equation}\\
(\breve \kappa_1-\kappa_1)-{\textstyle{1\over 6}}(\breve \kappa_3-\kappa_3)&=O(n^{-1}).\label{E:second equation}
\end{align}
Since $\eta T_1^2=(-1/\lambda^{11})\lambda^{1r}\lambda^{1s}l_rl_s=\tau^{rs}l_rl_s$, \eqref{E:first equation} yields
\begin{equation}
\eta^{1/2}[(\breve\xi^{rst}-\xi^{rst})l_{rs}l_t-(\breve\xi^{rs}-
\xi^{rs})l_rl_s-\{(\breve\xi^{tuv}-\xi^{tuv})\lambda_{tu,v}+
(\breve\xi^{tu}-\xi^{tu})\lambda_{tu}\}\tau^{rs}l_rl_s]=O_p(n^{-1}).\label{E:third equation}
\end{equation}
The quantity $\eta^{1/2}\{(\breve\xi^{rst}-\xi^{rst})l_{rs}l_t-(\breve\xi^{tuv}-\xi^{tuv})\lambda_{tu,v}\}$ in \eqref{E:third equation} is reduced to order $O_p(n^{-1})$ if \eqref{E:first condition} holds.

The remaining term $\eta^{1/2}\{(\breve\xi^{rs}-\xi^{rs})+(\breve\xi^{tu}-\xi^{tu})\lambda_{tu}\tau^{rs}\}l_rl_s$ in \eqref{E:third equation} is reduced to order $O_p(n^{-1})$ if
\eqref{E:second condition} holds. We show that \eqref{E:second equation} is satisfied when \eqref{E:first condition} and \eqref{E:second condition} hold.
Now \eqref{E:second equation} yields
\[
\eta^{1/2}\{(\breve\xi^{rst}-\xi^{rst})\lambda_{rs,t}+(\breve\xi^{rs}-\xi^{rs})\lambda_{rs}\}
+\eta^{3/2}\{(\breve\xi^{rs1}-\xi^{rs1})\lambda^{1t}\lambda_{rs,t}+\breve\xi^{11}-\xi^{11}\}=O(n^{-1}),
\]
and \eqref{E:first condition} yields $\breve\xi^{rs1}=\xi^{rs1}+O(n^{-5/2})$.  Under this condition, \eqref{E:second equation} reduces to
\[
\eta^{1/2}\{(\breve\xi^{rs}-\xi^{rs})\lambda_{rs}\}
+\eta^{3/2}(\breve\xi^{11}-\xi^{11})=O(n^{-1}).
\]
Since $\tau^{11}=-\lambda^{11}=\eta^{-1}$, \eqref{E:second condition} gives $\breve\xi^{11}-\xi^{11}=\eta^{-1}(\breve\xi^{rs}-\xi^{rs})\lambda_{rs}+O(n^{-5/2})$, and hence, it follows that under \eqref{E:first condition} and \eqref{E:second condition}, \eqref{E:second equation} is satisfied.
\endex

Note that \eqref{E:first equation} and \eqref{E:second equation} together constitute necessary and sufficient conditions for the $p$-values to agree to order $O_p(n^{-1})$. The quantity on the left side of \eqref{E:first equation} is of the form $\eta^{1/2}(A^{rst}l_{rs}l_t-B^{rs}l_rl_s)$, where
\[
A^{rst}=\breve\xi^{rst}-\xi^{rst}, \quad B^{rs}= (\breve\xi^{tuv}-\xi^{tuv})\lambda_{tu,v}\tau^{rs}+\breve\xi^{rs}-\xi^{rs}+(\breve\xi^{tu}-\xi^{tu})\lambda_{tu}\tau^{rs},
\]
so a necessary condition for agreement in general of $p$-values to order $O_p(n^{-1})$ is that $A^{rst}$ and $B^{rs}$ both be of order $O(n^{-5/2})$. The condition that $A^{rst}$ is of order $O(n^{-5/2})$ is the same as \eqref{E:first condition} and, in light of this condition, that $B^{rs}$ be of order $O(n^{-5/2})$ is equivalent to \eqref{E:second condition}. Thus, \eqref{E:first condition} and \eqref{E:second condition} are necessary for agreement of $p$-values to order $O_p(n^{-1})$. Of course, it is possible that the $p$-values from two test statistics $\breve T(\psi)$ and $T(\psi)$ fail to agree to order $O_p(n^{-1})$ for arbitrary models, yet they do agree for some specific model owing to particular features of the model. This situation could be revealed by verifying conditions \eqref{E:first condition} and \eqref{E:second condition} for the specific model. 

\section{Examples}

To illustrate the results of the previous sections, we consider eight asymptotically standard normal pivots, in addition to the signed root likelihood ratio statistic $R(\psi)$.

Consider four pivots that involve observed information.
For $R(\psi)$, we have $\xi^{rst}_R=\lambda^{1r}\lambda^{st}+{\textstyle{1\over 2}}\lambda^{1r}\tau^{st}$ and $\xi^{rs}_R={\textstyle{1\over 2}}\lambda^{1t}\lambda^{ru}\nu^{sv}\lambda_{tuv}+{\textstyle{1\over 6}}\lambda^{1t}\tau^{ru}\tau^{sv}\lambda_{tuv}$, and hence, $\xi^{rs}_R+\xi^{tu}_R\lambda_{tu}\tau^{rs}={\textstyle{1\over 2}}\lambda^{1t}\lambda^{ru}\nu^{sv}\lambda_{tuv}+{\textstyle{1\over 2}}\lambda^{1t}\nu^{uv}\lambda_{tuv}\tau^{rs}$.

Example 1. {\it Wald statistic with observed information.} For the Wald statistic defined by $T_{WO}(\psi)=(\hat\psi-\psi)\{-\hat M_{11}\}^{1/2}=(\hat\psi-\psi)\{-\hat L^{11}\}^{-1/2}$, we have $\xi^{rst}_{WO}=\xi^{rst}_R$ and  $\xi^{rs}_{WO}={\textstyle{1 \over 2}}\lambda^{1t}\lambda^{ru}\nu^{sv}\lambda_{tuv}$. Therefore, $\xi^{rs1}_{WO}={\textstyle {1\over 2}}\lambda^{1r}\lambda^{1s}$ and $\xi^{rs}_{WO}+\xi^{tu}_{WO}\lambda_{tu}\tau^{rs}=\xi^{rs}_R+\xi^{tu}_R\lambda_{tu}\tau^{rs}$.
We deduce that, to error of second order, $T_{WO}(\psi)$ is both stable in the sense discussed in Section 4 and produces the same $p$-values as $R(\psi)$.

Example 2. {\it Score statistic with observed information.} For the score statistic defined by $T_{SO}(\psi)=M_1(\psi)\{-\hat M_{11}\}^{-1/2}=L_1\{\tilde\theta(\psi)\}\{-\hat L^{11}\}^{1/2}$, we have $\xi^{rst}_{SO}=\xi^{rst}_R$ and  $\xi^{rs}_{SO}={\textstyle{1 \over 2}}\lambda^{1t}\lambda^{ru}\nu^{sv}\lambda_{tuv}+{\textstyle{1 \over 2}}\lambda^{1t}\tau^{ru}\tau^{sv}\lambda_{tuv}$. Thus, $\xi^{rs1}_{SO}={\textstyle {1\over 2}}\lambda^{1r}\lambda^{1s}$ and $\xi^{rs}_{SO}+\xi^{tu}_{SO}\lambda_{tu}\tau^{rs}=\xi^{rs}_R+\xi^{tu}_R\lambda_{tu}\tau^{rs}$.
It follows that, to error of second order, $T_{WO}(\psi)$ is also stable and again produces the same $p$-values as $R(\psi)$.

The following two asymptotically standard normal pivots are not standard components of likelihood-based inference. They involve pivots constructed by evaluating the observed information at the constrained maximum likelihood, rather than the global maximum likelihood estimator as in Examples 1 and 2.  Their use can be more cumbersome; they are included here to demonstrate the theoretical results.

Example 3. {\it Wald statistic with observed information evaluated at the constrained maximum likelihood estimator.} For the pivot $T_{WOC}(\psi)=(\hat\psi-\psi)[-M_{11}\{\tilde\theta(\psi)\}]^{1/2}=(\hat\psi-\psi)[-L^{11}\{\tilde\theta(\psi)\}]^{-1/2}$, we have $\xi^{rst}_{WOC}=\xi^{rst}_R$ and  $\xi^{rs}_{WOC}={\textstyle{1 \over 2}}\lambda^{1t}\lambda^{ru}\nu^{sv}\lambda_{tuv}+{\textstyle{1 \over 2}}\lambda^{1t}\tau^{ru}\tau^{sv}\lambda_{tuv}=\xi^{rs}_{SO}$. Hence, $\xi^{rs1}_{WOC}={\textstyle {1\over 2}}\lambda^{1r}\lambda^{1s}$ and $\xi^{rs}_{WOC}+\xi^{tu}_{WOC}\lambda_{tu}\tau^{rs}=\xi^{rs}_R+\xi^{tu}_R\lambda_{tu}\tau^{rs}$.
Thus $T_{WOC}(\psi)=T_{SO}(\psi)+O_p(n^{-1})$. To error of second order, $T_{WOC}(\psi)$ is stable and produces the same $p$-values as $R(\psi)$.

Example 4. {\it Score statistic with observed information evaluated at the constrained maximum likelihood estimator.} For $T_{SOC}(\psi)=M_1(\psi)[-M_{11}\{\tilde\theta(\psi)\}]^{-1/2}=L_1\{\tilde\theta(\psi)\}[-L^{11}\{\tilde\theta(\psi)\}]^{1/2}$, the corresponding score statistic, we have $\xi^{rst}_{SOC}=\xi^{rst}_R$ and  $\xi^{rs}_{SOC}={\textstyle{1 \over 2}}\lambda^{1t}\lambda^{ru}\nu^{sv}\lambda_{tuv}=\xi^{rs}_{WO}$. Thus, $\xi^{rs1}_{SOC}={\textstyle {1\over 2}}\lambda^{1r}\lambda^{1s}$ and $\xi^{rs}_{SOC}+\xi^{tu}_{SOC}\lambda_{tu}\tau^{rs}=\xi^{rs}_R+\xi^{tu}_R\lambda_{tu}\tau^{rs}$. As in the previous example, $T_{SOC}(\psi)=T_{WO}(\psi)+O_p(n^{-1})$. To error of second order, $T_{WOC}(\psi)$ is stable and produces the same $p$-values as $R(\psi)$.

We consider pivots corresponding to Examples 1-4 above, but based on expected, rather than observed, information.

Example 5. {\it Wald statistic with expected information.} For the version of the Wald statistic defined by $T_{WE}(\psi)=(\hat\psi-\psi)\{-\hat\lambda^{11}\}^{-1/2}$, we have $\xi^{rst}_{WE}=\lambda^{r1}\lambda^{st}$ and  $\xi^{rs}_{WE}={\textstyle{1 \over 2}}\lambda^{1t}\lambda^{ru}\nu^{sv}\lambda_{tuv}+{\textstyle{1 \over 2}}\lambda^{1t}\tau^{ru}\lambda^{sv}\lambda_{tu,v}$. Then, $\xi^{rs1}_{WE}=\lambda^{1r}\lambda^{1s}$ and $\xi^{rs}_{WE}+\xi^{tu}_{WE}\lambda_{tu}\tau^{rs}=\xi^{rs}_R+\xi^{tu}_R\lambda_{tu}\tau^{rs}+{\textstyle{1 \over 2}}\lambda^{1t}\tau^{ru}\lambda^{sv}\lambda_{tu,v}+{\textstyle{1 \over 2}}\lambda^{1t}\tau^{uv}\lambda_{tu,v}\tau^{rs}$.

Example 6. {\it Wald statistic with expected information evaluated at the constrained maximum likelihood estimator.} For the pivot described in Example 5, but with the expected information evaluated at the constrained maximum likelihood estimator, $T_{WEC}(\psi)=(\hat\psi-\psi)[-\lambda^{11}\{\tilde\theta(\psi)\}]^{-1/2}$, we have $\xi^{rst}_{WEC}=\xi^{rst}_{WE}$ and  $\xi^{rs}_{WEC}={\textstyle{1 \over 2}}\lambda^{1t}\lambda^{ru}\nu^{sv}\lambda_{tuv}+{\textstyle{1 \over 2}}\lambda^{1t}\lambda^{ru}\nu^{sv}\lambda_{tuv}+{\textstyle{1 \over 2}}\lambda^{1t}\tau^{ru}\tau^{sv}\lambda_{tu,v}$. Then, $\xi^{rs1}_{WEC}=\lambda^{1r}\lambda^{1s}$ and $\xi^{rs}_{WEC}+\xi^{tu}_{WEC}\lambda_{tu}\tau^{rs}=\xi^{rs}_{WE}+\xi^{tu}_{WE}\lambda_{tu}\tau^{rs}$.

Neither $T_{WE}(\psi)$ nor $T_{WEC}(\psi)$ generally satisfy the above sufficient condition for stability to error of order $O(n^{-1})$ and, of course, they do not generally provide $p$-values that agree with those from $R(\psi)$ to error of order $O_p(n^{-1})$. However, the $p$-values calculated from $T_{WE}(\psi)$ agree with those from  $T_{WEC}(\psi)$ to error of order $O_p(n^{-1})$.

Example 7. {\it Score statistic with expected information.} For the version of the score statistic defined by $T_{SE}(\psi)=M_1(\psi)\{-\hat\lambda^{11}\}^{1/2}=\nobreak  L_1\{\tilde\theta(\psi)\}\nobreak\{-\hat\lambda^{11}\}^{1/2}$, we have $\xi^{rst}_{SE}=\lambda^{r1}\nu^{st}$ and  $\xi^{rs}_{SE}={\textstyle{1 \over 2}}\lambda^{1t}\lambda^{ru}\nu^{sv}\lambda_{tuv}+{\textstyle{1 \over 2}}\lambda^{1t}\tau^{ru}\tau^{sv}\lambda_{tuv}-{\textstyle{1 \over 2}}\lambda^{1t}\tau^{ru}\lambda^{sv}\lambda_{tu,v}$. Therefore, $\xi^{rs1}_{SE}=0$ and $\xi^{rs}_{SE}+\xi^{tu}_{SE}\lambda_{tu}\tau^{rs}=\xi^{rs}_R+\xi^{tu}_R\lambda_{tu}\tau^{rs}-{\textstyle{1 \over 2}}\lambda^{1t}\tau^{ru}\lambda^{sv}\lambda_{tu,v}-{\textstyle{1 \over 2}}\lambda^{1t}\tau^{uv}\lambda_{tu,v}\tau^{rs}$.

Example 8. {\it Score statistic with expected information evaluated at the constrained maximum likelihood estimator.} Evaluating the expected information instead at the constrained maximum likelihood estimator, for $T_{SEC}(\psi)=M_1(\psi)[-\lambda^{11}\{\tilde\theta(\psi)\}]^{1/2}=L_1\{\tilde\theta(\psi)\}[-\lambda^{11}\{\tilde\theta(\psi)\}]^{1/2}$, we have $\xi^{rst}_{SE}=\lambda^{r1}\nu^{st}$ and  $\xi^{rs}_{SEC}={\textstyle{1 \over 2}}\lambda^{1t}\lambda^{ru}\nu^{sv}\lambda_{tuv}-{\textstyle{1 \over 2}}\lambda^{1t}\tau^{ru}\nu^{sv}\lambda_{tu,v}$. Thus, $\xi^{rs1}_{SEC}=0$ and $\xi^{rs}_{SEC}+\xi^{tu}_{SEC}\lambda_{tu}\tau^{rs}=\xi^{rs}_{SE}+\xi^{tu}_{SE}\lambda_{tu}\tau^{rs}$.

Neither $T_{SE}(\psi)$ nor $T_{SEC}(\psi)$ generally satisfy the above sufficient condition for stability to error of order $O(n^{-1})$, and they do not generally provide $p$-values that agree with those from $R(\psi)$ to error of order $O_p(n^{-1})$. However, the $p$-values calculated from $T_{SE}(\psi)$ agree with those from  $T_{SEC}(\psi)$ to error of order $O_p(n^{-1})$, although they do not generally agree with those from $T_{WE}(\psi)$ and $T_{WEC}(\psi)$ to error of order $O_p(n^{-1})$.

Construction of the asymptotically normal pivot for inference on the interest parameter $\psi$ in the presence of a nuisance parameter using observed information is therefore key to ensuring that $p$-values calculated from the marginal distribution of the pivot, as might be approximated in generality by parametric bootstrapping, automatically respect, to second-order, the conditioning on ancillary statistics required for inferential correctness. The importance of using observed information instead of expected information for approximate conditional inference is, of course, well known, having been argued by Efron and Hinkley (1978), who were partly inspired by the discussion given by Pierce (1975) to the paper by Efron (1975) on the geometry of exponential families. Our analysis gives a very direct operational interpretation, in terms of the $p$-values derived from the marginal sampling distributions of commonly used pivots.

Further discrimination between pivots may be based on the requirement of parameterisation invariance, that inferential conclusions should not depend on the parameterisation: see, for instance, Pace and Salvan (1997, Section 2.11).  Requirement of invariance of the inference under reparameterisations which are (Barndorff-Nielsen and Cox (1994, Section 1.5)) interest-respecting would exclude use of Wald statistics: see, for instance, McCullagh (1987, Section 7.4).

\section{Extension to adjusted profile likelihood}

The general form of the asymptotically normal test statistic that we have considered, where the statistic is expressible as $T(\psi)=\eta^{1/2}(T_1+T_2)+O_p(n^{-1})$, where $T_1=-\lambda^{1r}l_r$ and $T_2$ is of the form $T_2=\xi^{rst}l_{rs}l_t-\xi^{rs}l_rl_s$, with $\xi^{rst}$ and $\xi^{rs}$ assumed to be of order $O(n^{-2})$, covers important special cases which are commonly applied. It does not, however, include asymptotically standard normal pivots based on adjusted forms of profile likelihood. Fortunately, only a simple change to the analysis is necessary is accommodate pivots based on adjusted likelihoods. The criteria for second-order stability and equivalence of $p$-values are unchanged since, to the order being considered, the version of the pivot based on the adjusted profile likelihood is obtained by a constant, additive adjustment of that based on the unadjusted profile likelihood.

There have been many suggestions to replace the usual profile likelihood function $M(\psi)$ by an adjusted version $\bar M(\psi)=M(\psi)+B(\psi)$, where $B(\psi)$ is an adjustment function which is a function of $Y$ and $\psi$ only, whose derivatives with respect to $\psi$ are of order $O_p(1)$.
The likelihood ratio statistic based on the adjusted profile likelihood is $\bar W(\psi)=2\{\bar M(\bar\psi)-\bar M(\psi)\}$, where $\bar\psi$ is the point at which $\bar M(\psi)$ is maximized.
The signed root of the likelihood ratio statistic based on the adjusted profile likelihood is $\bar R(\psi)={\rm sgn}(\bar\psi-\psi)\{\bar W(\psi)\}^{1/2}$.

Following our previous notation, we write $B_1(\psi)=\partial B(\psi)/\partial\psi$, $B_{11}(\psi)=\partial^2 B(\psi)/\partial\psi^2$, etc.
Let $\beta_1=E\{B_1(\psi)\}$, $\beta_{11}=E(B_{11})$, etc.; these quantities are assumed to be of order $O(1)$.
Further, let $b_1=B_1(\psi)-\beta_1$, $b_{11}=B_{11}(\psi)-\beta_{11}$, etc., with these quantities assumed to be of order $O_p(n^{-1/2})$.
Assume also that the joint cumulants of $nb_1$, $nb_{11}$, $l_r$, $l_{rs}$, etc.\ are of order $O(n)$.

In many instances, a specific adjustment function $B(\psi)$ has been proposed to take into account the effect of nuisance parameters for inference about $\psi$, notably the modified profile likelihood of Barndorff-Nielsen (1983) and the adjusted profile likelihood of Cox and Reid (1987). Other adjustments with the same structure as described above are detailed by Skovgaard (1996), Severini (1998), DiCiccio and Martin (1993), and Barndorff-Nielsen and Chamberlin (1994). These adjustment functions have the effect of reducing the mean of the profile score from order $O(1)$ to order $O(n^{-1})$: see, for instance, DiCiccio et al.\ (1996).
The adjustment functions have $\beta_1=\rho+O(n^{-1})$, where $\rho=-\eta\lambda^{1r}\nu^{st}({\textstyle{1\over 2}}\lambda_{rst}+\lambda_{rs,t})$. Since, in general, $E\{M_1(\psi)\}=-\rho+O(n^{-1})$, it follows that $E\{\bar M_1(\psi)\}=O(n^{-1})$: see McCullagh and Tibshirani (1990), DiCiccio et al.\ (1996).

Another version of the adjustment function that derives from Bayesian inference based on a prior density $\pi(\theta)$ is
\[
B(\psi)=-{1\over 2}\log\biggl(\frac{\det[-L_{ab}\{\tilde\theta(\psi)\}]}{\det\{-L_{ab}(\hat\theta)\}}\biggr)+\log\biggl[\frac{\pi\{\tilde\theta(\psi)\}}{\pi(\hat\theta)}\biggr],
\]
where $a,b=2,\ldots,d$. Here $\{L_{ab}(\theta)\}$ is the $(d-1)\times(d-1)$ submatrix of $\{L_{rs}(\theta)\}$ corresponding to the nuisance parameters.
This adjustment function arises from the Laplace approximation to $\pi_{\psi|Y}(\psi)$, the posterior marginal density function for $\psi$, developed by Tierney and Kadane (1986), who showed that $\pi_{\psi|Y}(\psi)=c\bar M(\psi)\{1+O(n^{-3/2})\}$, for values of $\psi$ such that $\psi-\hat\psi$ is of order $O(n^{-1/2})$.
In this case, $\bar W(\psi)$ corresponds to the posterior ratio statistic to error of order $O_p(n^{-3/2})$, and $\beta_1=\eta\lambda^{1r}({\textstyle{1 \over 2}}\nu^{st}\lambda_{rst}-\pi_r/\pi)$: see DiCiccio and Stern (1994a).  Firth (1993) developed particular adjustment functions motivated by the specific aim that $\bar\psi$ be unbiased to error of order $O(n^{-3/2})$.

For a general adjustment function $B(\psi)$, DiCiccio and Stern (1994a) showed that $\bar R(\psi)=\eta^{1/2}\{\bar R_1+\bar R_2+O_p(n^{-3/2})\}$, where $\bar R_1=R_1=-\lambda^{1r}l_r$ and $\bar R_2=R_2-\lambda^{11}\beta_1$; in particular, $\bar R(\psi)=R(\psi)+\eta^{-1/2}\beta_1+O_p(n^{-1})$.

Pierce and Bellio (2006), considering the adjustment functions related to modified profile likelihood and Bayesian inference, also observed that, to error of order $O_p(n^{-1})$, $\bar R(\psi)$ differs from $R(\psi)$ by only a constant, although they did not detail the associated formulae involving $\beta_1$.
Having made this observation, Pierce and Bellio (2006) conclude that, to error of order $O_p(n^{-1})$, both $\bar R(\psi)$ and $R(\psi)$ induce the same orderings of datasets for evidence against the null hypothesis, and they conclude that, to this order of error, ideal frequentist $p$-values can be based on the distribution of $R(\psi)$.

We generalize our preceding results by considering hypothesis testing for $\psi$ based on a test statistic $\bar T(\psi)=\eta^{1/2}(\bar T_1+\bar T_2)+O_p(n^{-1})$ where, as before, $\bar T_1=T_1=-\lambda^{1r}l_r$, and $\bar T_2$ is assumed to be of the form $\bar T_2=\xi^{rst}l_{rs}l_t-\xi^{rs}l_rl_s+\varsigma=T_2+\varsigma$, with $\xi^{rst}$ and $\xi^{rs}$ of order $O(n^{-2})$ and the constant $\varsigma$ assumed to be of order $O(n^{-1})$. Therefore, $\bar T(\psi)=T(\psi)+\eta^{1/2} \varsigma+O(n^{-1})$. We provide illustrations which demonstrate how statistics constructed from adjusted profile likelihood may be expressed in this form.

Since $\bar T(\psi)$ only differs, to the second-order being considered, from $T(\psi)$ by a constant, the condition for $\bar T(\psi)$ to be stable to error of order $O(n^{-1})$ is the same as the condition for $T(\psi)$, namely $\xi^{rs1}={\textstyle{1\over 2}}\lambda^{1r}\lambda^{1s}$.

The first three cumulants of $\bar T(\psi)=T(\psi)+\eta^{1/2} \varsigma+O(n^{-1})$ are $\bar\kappa_1=\kappa_1+\eta^{1/2} \varsigma+O(n^{-1})$, $\bar\kappa_2=\kappa_2+O(n^{-1})$, $\bar\kappa_3=\kappa_3+O(n^{-1})$, where $\kappa_1$, $\kappa_2$, and $\kappa_3$ are as described before for $T(\psi)$, and the fourth- and higher-order cumulants of $\bar T(\psi)$ are of order $O(n^{-1})$, or smaller.

Consider two versions of $\bar T(\psi)$, say $T(\psi)+\eta^{1/2} \varsigma+O(n^{-1})$ and $\breve T(\psi)+\eta^{1/2} \breve\varsigma+O(n^{-1})$.
The preceding Cornish-Fisher argument for comparing $p$-values shows that the $p$-values from the two test statistics differ by order $O_p(n^{-1})$ provided
\[
\{\eta^{1/2}(\breve T_2+ \breve\varsigma-T_2- \varsigma)-{\textstyle{1\over 6}}(\breve \kappa_3-\kappa_3)\eta T_1^2\}
-\{(\breve \kappa_1+\eta^{1/2} \breve\varsigma-\kappa_1-\eta^{1/2} \varsigma)-{\textstyle{1\over 6}}(\breve \kappa_3-\kappa_3)\}=O_p(n^{-1}).
\]
The crucial point is that the terms involving $\varsigma$ and $\breve\varsigma$ cancel from the left side of this expression, irrespective of their values, so \eqref{E:first condition} and \eqref{E:second condition} continue to specify necessary and sufficient conditions for the two test statistics to yield $p$-values that differ by order $O_p(n^{-1})$.

Example 9. {\it Signed root likelihood ratio statistic constructed from adjusted profile likelihood.} For the signed root likelihood ratio statistic constructed from the adjusted profile likelihood, $\bar R(\psi)$, standard calculations show that $\xi^{rst}_{\bar R}=\xi^{rst}_R$, $\xi^{rs}_{\bar R}=\xi^{rs}_R$, $\varsigma_{\bar R}=\eta^{-1}\beta_1$. It follows that, to error of order $O_p(n^{-1})$, $\bar R(\psi)$ and $R(\psi)$ produce the same $p$-values, as noted by Pierce and Bellio (2006).

Example 10. {\it Wald statistic with observed information constructed from adjusted profile likelihood.} For the pivot $T_{AWO}(\psi)=(\bar\psi-\psi)\{-\bar M_{11}(\bar\psi)\}^{1/2}$, we have $\xi^{rst}_{AWO}=\xi^{rst}_{WO}=\xi^{rst}_R$,  $\xi^{rs}_{AWO}=\xi^{rs}_{WO}$, and $\varsigma_{AWO}=\eta^{-1}\beta_1$. Then, since to error of order $O_p(n^{-1})$, $T_{WO}(\psi)$ and $R(\psi)$ produce the same $p$-values, it follows that $T_{AWO}(\psi)$ and $R(\psi)$ produce the same $p$-values to that order of error.

Example 11. {\it Score statistic with observed information constructed from adjusted profile likelihood.} For the statistic $T_{ASO}(\psi)=\bar M_1(\psi)\{-\bar M_{11}(\bar\psi)\}^{1/2}$, we have $\xi^{rst}_{ASO}=\xi^{rst}_{SO}=\xi^{rst}_R$,  $\xi^{rs}_{ASO}=\xi^{rs}_{SO}$, and $\varsigma_{ASO}=\eta^{-1}\beta_1$. Since, to error of order $O_p(n^{-1})$, $T_{SO}(\psi)$ and $R(\psi)$ produce the same $p$-values, it follows that $T_{ASO}(\psi)$ and $R(\psi)$ produce the same $p$-values to that order of error.

The interesting feature here is that although $\bar R(\psi)$, $T_{AWO}(\psi)$, and $T_{ASO}(\psi)$ differ from one another by non-constant terms of order $O_p(n^{-1/2})$ in general, they all produce the same $p$-values to error of order $O_p(n^{-1})$.

\section{Vector-valued interest parameter}

Consider again the partition $\theta=(\psi,\phi)$, but now allow for the possibility that the interest parameter $\psi$ is vector-valued, having dimension $q$.
The likelihood ratio statistic $W(\psi)$ is routinely used for hypothesis testing about $\psi$.
The asymptotic distribution of $W(\psi)$ is chi-squared with $q$ degrees of freedom. Indeed, for regular problems, the $\chi^2_q$-approximation to the distribution of $W(\psi)$ has error of order $O(n^{-1})$, and moreover, the mean of $W(\psi)$ has the expansion $E\{W(\psi)\}=q(1+n^{-1}\omega)+O(n^{-2})$, where $\omega \equiv \omega(\theta)$ is of order $O(1)$.
Lawley (1956), Barndorff-Nielsen and Cox (1984), and Bickel and Ghosh (1990) showed that $W(\psi)$ is distributed as $(1+n^{-1}\omega)\chi^2_q$ to error of order $O(n^{-2})$: the Bartlett-corrected statistic $W(\psi)/(1+n^{-1} \omega)$ is distributed as $\chi^2_q$ to error of order $O(n^{-2})$. Further, $W(\psi)$ is stable.
\begin{theorem}
The unconditional and conditional distributions of $W(\psi)$ agree to error of order $O(n^{-3/2})$, given the ancillary statistic $A$.
\end{theorem}
\emph{Proof.}
By applying identical arguments to the conditional distribution of $Y$ given $A$, we have that $\mathring E\{W(\psi)\}=q(1+n^{-1}\mathring \omega)+O(n^{-2})$, where $\mathring \omega$ is of order $O(1)$ given $A$, and that $W(\psi)$ is conditionally distributed as $(1+n^{-1}\mathring\omega)\chi^2_q$ to error of order $O(n^{-2})$ given $A$.
Barndorff-Nielsen and Cox (1984) showed that $\mathring\omega=\omega+O_p(n^{-1/2})$, and hence it follows that $W(\psi)$ is stable to error of order $O(n^{-3/2})$. Extending the arguments of McCullagh (1987, Section 8.4) to the nuisance parameter case, $\mathring\omega=\omega+O_p(n^{-1/2})$ continues to hold provided the conditioning statistic $A$ is a second-order local ancillary statistic.
\endex

Inference based on an approximation to the marginal distribution of $W(\psi)$ accurate to error of order $O(n^{-3/2})$ therefore automatically respects conditioning on the ancillary statistic to that same order.

Bickel and Ghosh (1990) explicitly recommended that the Bartlett adjustment factor  $(1+n^{-1}\omega)$ be estimated by simulation; this may be done by either fixing $\theta=\hat\theta$ or $\theta=\tilde\theta$, so that inference is based on a $\chi^2_q$ approximation to the sampling distribution of, say, $W(\psi)/\{1+n^{-1} \omega(\tilde \theta)\}$. Alternatively, the entire distribution of $W(\psi)$ may be approximated by simulation at either of these parameter values: such an approximation is, however, likely to be computationally more expensive than estimation of just the Bartlett adjustment factor. In view of the stability result above, these inference procedures not only provide $p$-values that are uniformly distributed to error of order $O_p(n^{-3/2})$ (actually, the error is of order $O_p(n^{-2})$ - see Barndorff-Nielsen and Hall (1988)), but these $p$-values are uniformly distributed conditionally to the same order of error.

DiCiccio and Stern (1994b) demonstrated the efficacy of Bartlett correction for likelihood ratio statistics based on adjusted profile likelihoods. They showed that $E\{\bar W(\psi)\}=q(1+n^{-1}\bar\omega)+O(n^{-2})$ and that $\bar W(\psi)$ is distributed as $(1+n^{-1}\bar\omega)\chi^2_q$ to error of order $O(n^{-2})$. Moreover, their calculations can be applied to the conditional distribution of $Y$ given $A$ to show that these results also hold conditionally, as for $W(\psi)$.

\begin{theorem}
The unconditional and conditional distributions of $\bar W(\psi)$ agree to order $O(n^{-3/2})$, given the ancillary statistic $A$.
\end{theorem}
\emph{Proof.}
See Appendix.
\endex

The operational consequences of this stability result are again straightforward.  Similar stability results hold for other test statistics that are asymptotically distributed as $\chi^2_q$, such as $(\bar \psi^a-\psi^a)(\bar \psi^b-\psi^b)\bar S_{ab}$ and $\bar M_a(\psi)\bar M_b(\psi)\bar S^{ab}$, where
$\bar S_{ab}=-\bar M_{ab}(\bar \psi)$ and $(\bar S^{ab})$ is the $q \times q$ matrix inverse of $(\bar S_{ab})$. The marginal distribution function of such a statistic $X$ typically has the expansion
\[
Pr(X \leq x)=Pr(\chi^2_q \leq x) +\sum_{j=0}^k \alpha_j Pr(\chi^2_{q+2j} \leq x)+O(n^{-3/2}),
\]
where the $\alpha_j$ are functions of the $\lambda$'s and $\beta$'s and typically $k=3$; see, for example, Harris (1985) and Cordeiro and Ferrari (1991). The same manipulations of likelihood quantities that produce the approximation to the marginal distribution of $X$ can be applied to conditional likelihood quantities to yield the expansion
\[
Pr(X \leq x \mid A)=Pr(\chi^2_q \leq x) +\sum_{j=0}^k \mathring \alpha_j Pr(\chi^2_{q+2j} \leq x)+O_p(n^{-3/2}),
\]
where the $\mathring \alpha_j$ are functions of the $\mathring \lambda$'s and $\mathring \beta$'s. The preceding calculations that demonstrate the stability of $\bar W(\psi)$ can also be used to show that $\mathring \alpha_j=\alpha_j+O_p(n^{-3/2}),$ and it follows that $X$ is stable to error of order $O(n^{-3/2})$.

\section{Discussion}

Focus here has been on inference on an interest parameter in the presence of a nuisance parameter in ancillary statistic models. We have shown that commonly used, asymptotically standard normal, likelihood-based pivots, including the signed root statistic $R(\psi)$, are second-order stable. When applied with such a pivot, procedures such as the parametric bootstrap, which approximate the marginal distribution of the pivot to second-order, achieve the same order of accuracy, $O(n^{-1})$, in approximation of the relevant exact conditional inference. Our motivation for the analysis here is as a preliminary to full evaluation of the properties of such parametric bootstrap procedures as an alternative to more awkward analytic approaches to approximation of exact conditional inference. In this regard, of importance for future investigation is analysis of large deviation properties of procedures based on marginal simulation of a likelihood-based pivot. Analytic procedures, such as normal approximation to $R^*(\psi)$, or the approximation of Skovgaard (1996), confer large deviation protection, typically providing accurate approximation of the conditional distribution of the associated pivot far into its tails. The requirement of such large deviation behaviour may be judged an important discriminant between competing methodologies. Discussion of this and related issues is currently in preparation in DiCiccio, Kuffner and Young (2014).

Pivots stable to third-order do, of course, exist: $R^*(\psi)$ is distributed as standard normal to third-order, conditionally on the ancillary statistic, and hence unconditionally as well. Second-order approximation to an exact conditional inference through the bootstrap is seen (see, for example, DiCiccio and Young (2010), Young and Smith (2005, Chapter 10)) to give good results in practice in ancillary statistic settings. Basing inference on a pivot stable to third-order seems unwarranted. In addition, ancillary statistics are typically not unique and (see, for instance, McCullagh (1992)), different conditional inferences typically only agree to second-order, so it can be argued that third-order approximation to an exact conditional inference is, in itself, unwarranted. By our analysis, inference based on second-order (or higher-order) approximation of the marginal distribution of a pivot stable to second-order approximates {\it any} conditional inference to $O(n^{-1})$.

Our study of uniqueness of $p$-values yielded simple conditions under which $p$-values derived from different asymptotically standard normal pivots agree to order $O_{p}(n^{-1})$. In cases we have considered where the conditions fail to be satisfied, a more detailed analysis shows that $p$-values agree only to an actual order $O_p(n^{-1/2})$.

\section*{Appendix}
\subsection*{Proof of Lemma 2}

The unconditional variance of $T(\psi)$ is
\begin{align*}
{\rm var}\{T(\psi)\}&=E[\{T(\psi)\}^2]-[E\{T(\psi)\}]^2=E[\{T(\psi)\}^2]+O(n^{-1})\\
&=\eta E\{T_1^2+2T_1T_2+O_p(n^{-2})\}+O(n^{-1})\\
&=\eta E\{\lambda^{1r}\lambda^{1s}l_rl_s-2\lambda^{1r}\xi^{stu}l_rl_{st}l_u+2\lambda^{1r}\xi^{st}l_rl_sl_t+O_p(n^{-2})\}+O(n^{-1})\\
&=-\eta\{\lambda^{1r}\lambda^{1s}\lambda_{rs}+O(n^{-2})\}+O(n^{-1})\\
&=1+O(n^{-1}).
\end{align*}
Correspondingly, the conditional variance of $T(\psi)$ is
\begin{align*}
\mathring{{\rm var}}\{T(\psi)\}&=\mathring E[\{T(\psi)\}^2]-[\mathring E\{T(\psi)\}]^2=\mathring E[\{T(\psi)\}^2]+O_p(n^{-1})\\
&=\eta \mathring E\{T_1^2+2T_1T_2+O_p(n^{-2})\}+O_p(n^{-1})\\
&=\eta \mathring E\{\lambda^{1r}\lambda^{1s}\mathring l_r \mathring l_s-2\lambda^{1r}\xi^{stu}\mathring l_r(\mathring l_{st}+\mathring\Delta_{st})\mathring l_u+2\lambda^{1r}\xi^{st}\mathring l_r\mathring l_s\mathring l_t+O_p(n^{-2})\}+O_p(n^{-1})\\
&=-\eta\{\lambda^{1r}\lambda^{1s}\mathring\lambda_{rs}-2\lambda^{1r}\xi^{stu}\mathring\lambda_{ru}\mathring\Delta_{st}+O_p(n^{-2})\}+O_p(n^{-1})\\
&=-\eta\{\lambda^{1r}\lambda^{1s}(\lambda_{rs}+\mathring\Delta_{rs})-2\lambda^{1r}\xi^{stu}\lambda_{ru}\mathring\Delta_{st}\}+O_p(n^{-1})\\
&=1-\eta(\lambda^{1r}\lambda^{1s}\mathring\Delta_{rs}-2\xi^{st1}\mathring\Delta_{st})+O_p(n^{-1})\\
&=1-\eta\{(\lambda^{1r}\lambda^{1s}-2\xi^{rs1})\mathring\Delta_{rs}\}+O_p(n^{-1}).\\
\end{align*}
It follows that $\mathring{{\rm var}}\{T(\psi)\}={\rm var}\{T(\psi)\}+O_p(n^{-1})$ provided $\xi^{rs1}={\textstyle{1\over 2}}\lambda^{1r}\lambda^{1s}$.
\endex
\par

\subsection*{Proof of Lemma 3}

The unconditional skewness of $T(\psi)$ is
\begin{align*}
{\rm skew}\{T(\psi)\}&=E([T(\psi)-E\{T(\psi)\}]^3)=E[\{T(\psi)\}^3]-3E[\{T(\psi)\}^2]E\{T(\psi)\}+O(n^{-1})\\
&=\eta^{3/2}[E\{(T_1+T_2)^3\}-3E\{(T_1+T_2)^2\}E(T_1+T_2)]+O(n^{-1})\\
&=\eta^{3/2}[E\{T_1^3+3T_1^2T_2+O_p(n^{-5/2})\}-3E\{T_1^2+O_p(n^{-3/2})\}E(T_2)]+O(n^{-1})\\
&=\eta^{3/2}[E\{-\lambda^{1r}\lambda^{1s}\lambda^{1t}l_rl_sl_t
+3\lambda^{1r}\lambda^{1s}(\xi^{tuv}l_{tu}l_v-\xi^{tu}l_tl_u)l_rl_s+O_p(n^{-5/2})\}\\
&\hspace{20pt}-3E\{\lambda^{1r}\lambda^{1s}l_rl_s+O_p(n^{-3/2})\}\{\xi^{rst}\lambda_{rs,t}+\xi^{rs}\lambda_{rs}+O_p(n^{-3/2})\}]+O(n^{-1})\\
&=\eta^{3/2}\{E(-\lambda^{1r}\lambda^{1s}\lambda^{1t}l_rl_sl_t+3\lambda^{1r}\lambda^{1s}\xi^{tuv}l_rl_sl_{tu}l_v
-3\lambda^{1r}\lambda^{1s}\xi^{tu}l_rl_sl_tl_u\\
&\hspace{20pt}-3\lambda^{1r}\lambda^{1s}\xi^{tuv}l_rl_s\lambda_{tu,v}-3\lambda^{1r}
\lambda^{1s}\xi^{tu}l_rl_s\lambda_{tu})\}+O(n^{-1}).
\end{align*}
To continue the calculation, we make use of the following identities:
\begin{align*}
-E(l_rl_sl_t)&=\lambda_{rs,t}+\lambda_{rt,s}+\lambda_{st,r}+\lambda_{rst},\\
E(l_rl_sl_{tu}l_v)&=-\lambda_{rs}\lambda_{tu,v}-\lambda_{rv}\lambda_{tu,s}-\lambda_{sv}\lambda_{tu,r}+O(n^{3/2}),\\
E(l_rl_sl_tl_u)&=\lambda_{rs}\lambda_{tu}+\lambda_{rt}\lambda_{su}+\lambda_{ru}\lambda_{st}+O(n^{3/2}).
\end{align*}
By using these identities, we obtain
\begin{align*}
{\rm skew}\{T(\psi)\}&=\eta^{3/2}(3\lambda^{1r}\lambda^{1s}\lambda^{1t}\lambda_{rs,t}+\lambda^{1r}\lambda^{1s}\lambda^{1t}\lambda_{rst}\\
&\hspace{20pt}-3\lambda^{11}\xi^{tuv}\lambda_{tu,v}-3\lambda^{1s}\xi^{tu1}\lambda_{tu,s}-3\lambda^{1r}\xi^{tu1}\lambda_{tu,r}\\
&\hspace{40pt}-3\lambda^{11}\xi^{tu}\lambda_{tu}-3\xi^{11}-3\xi^{11}\\
&\hspace{60pt}+3\lambda^{11}\xi^{tuv}\lambda_{tu,v}+3\lambda^{11}\xi^{tu}\lambda_{tu})+O(n^{-1})\\
&=\eta^{3/2}(\lambda^{1r}\lambda^{1s}\lambda^{1t}\lambda_{rst}+3\lambda^{1r}\lambda^{1s}\lambda^{1t}\lambda_{rs,t}
-6\xi^{rs1}\lambda^{1t}\lambda_{rs,t}-6\xi^{11})+O(n^{-1}).
\end{align*}
Similar reasoning shows that the conditional skewness of $T(\psi)$ is
\begin{align*}
\mathring{{\rm skew}}\{T(\psi)\} &=\eta^{3/2}[\mathring E\{-\lambda^{1r}\lambda^{1s}\lambda^{1t}l_rl_sl_t+3\lambda^{1r}\lambda^{1s}\xi^{tuv}l_rl_sl_{tu}l_v -3\lambda^{1r}\lambda^{1s}\xi^{tu}l_rl_sl_tl_u\\
&\hspace{20pt}-3\lambda^{1r}\lambda^{1s}\xi^{tuv}l_rl_s\lambda_{tu,v}-3\lambda^{1r}\lambda^{1s}\xi^{tu}
l_rl_s\lambda_{tu}+O_p(n^{-5/2})\}]+O_p(n^{-1})\\
&=\eta^{3/2}[\mathring E\{-\lambda^{1r}\lambda^{1s}\lambda^{1t}\mathring l_r\mathring l_s\mathring l_t+3\lambda^{1r}\lambda^{1s}\xi^{tuv}\mathring l_r\mathring l_s(\mathring l_{tu}+\mathring \Delta_{tu})\mathring l_v -3\lambda^{1r}\lambda^{1s}\xi^{tu}\mathring l_r\mathring l_s\mathring l_t\mathring l_u\\
&\hspace{20pt}-3\lambda^{1r}\lambda^{1s}\xi^{tuv}\mathring l_r\mathring l_s\lambda_{tu,v}-3\lambda^{1r}\lambda^{1s}\xi^{tu}\mathring l_r\mathring l_s\lambda_{tu}\}]+O_p(n^{-1})\\
&=\eta^{3/2}\{\mathring E(-\lambda^{1r}\lambda^{1s}\lambda^{1t}\mathring l_r\mathring l_s\mathring l_t+3\lambda^{1r}\lambda^{1s}\xi^{tuv}\mathring l_r\mathring l_s\mathring l_{tu}\mathring l_v -3\lambda^{1r}\lambda^{1s}\xi^{tu}\mathring l_r\mathring l_s\mathring l_t\mathring l_u\\
&\hspace{20pt}-3\lambda^{1r}\lambda^{1s}\xi^{tuv}\mathring l_r\mathring l_s\lambda_{tu,v}-3\lambda^{1r}\lambda^{1s}\xi^{tu}\mathring l_r\mathring l_s\lambda_{tu})\}+O_p(n^{-1}).
\end{align*}
Now we use the following identities:
\begin{align*}
-\mathring E(\mathring l_r\mathring l_s\mathring l_t)&=\mathring \lambda_{rs,t}+\mathring \lambda_{rt,s}+\mathring \lambda_{st,r}+\mathring \lambda_{rst}\\
&=\lambda_{rs,t}+\lambda_{rt,s}+\lambda_{st,r}+\lambda_{rst} + O_p(n^{1/2})\\
&=\lambda_{r,s,t}+ O_p(n^{1/2})\\
&=-E(l_rl_sl_t)+ O_p(n^{1/2}),\\[10pt]
\mathring E(\mathring l_r\mathring l_s\mathring l_{tu}\mathring l_v)&=-\mathring \lambda_{rs}\mathring \lambda_{tu,v}-\mathring \lambda_{rv}\mathring \lambda_{tu,s}-\mathring \lambda_{sv}\mathring \lambda_{tu,r}+O_p(n^{3/2})\\
&=-\lambda_{rs}\lambda_{tu,v}-\lambda_{rv}\lambda_{tu,s}-\lambda_{sv}\lambda_{tu,r}+O_p(n^{3/2})\\
&=E(l_r l_s l_{tu} l_v)+O_p(n^{3/2}),\\[10pt]
\mathring E(\mathring l_r\mathring l_s\mathring l_t\mathring l_u)&=\mathring \lambda_{rs}\mathring \lambda_{tu}+\mathring \lambda_{rt}\mathring \lambda_{su}+\mathring \lambda_{ru}\mathring \lambda_{st}+O_p(n^{3/2})\\
&=\lambda_{rs}\lambda_{tu}+\lambda_{rt}\lambda_{su}+\lambda_{ru}\lambda_{st}+O_p(n^{3/2})\\
&= E( l_r l_s l_t l_u)+O_p(n^{3/2}).
\end{align*}
By using these identities in the preceding expression for $\mathring{{\rm skew}}\{T(\psi)\}$, it is apparent that $\mathring{{\rm skew}}\{T(\psi)\}={{\rm skew}}\{T(\psi)\}+O_p(n^{-1})$, and hence, the conditional third cumulant agrees with the unconditional one to error of order $O_p(n^{-1})$, as required.
\endex
\par

\subsection*{Proof of Theorem 5}

To establish the stability of $\bar W(\psi)$ to error of order $O(n^{-3/2})$, we need only show that $\mathring E\{\bar W(\psi)\}=E\{\bar W(\psi)\}+O_p(n^{-3/2})$. For full generality, the previous notation, which is applicable when $\psi$ is a scalar, must be extended. In the expressions that follow, it is assumed that
subscripts and superscripts $a, b, \ldots$ have the range $1, \ldots, q$, while $r, s, \ldots$ range over $1, \ldots, d$. Let $(\eta_{ab})$ be the $q \times q$ matrix inverse of $(-\lambda^{ab})$, let $\tau^{rs}=\eta_{ab}\lambda^{ar}\lambda^{bs}$, and let $\nu^{rs}=\lambda^{rs}+\tau^{rs}$. In addition, let $B_a(\psi)=\partial B(\psi)/\partial \psi^a, B_{ab}(\psi)=\partial^2 B(\psi)/\partial \psi^a \partial \psi^b,
\beta_a=E\{B_a(\psi)\}, \beta_{ab}=E\{B_{ab}(\psi)\}, b_a=B_a(\psi)-\beta_a, b_{ab}=B_{ab}(\psi)-\beta_{ab}$, and so forth. The constants $\beta_a, \beta_{ab}$ etc. are assumed to be of order $O(1)$ and the variables $b_a, b_{ab}$ etc. are assumed to be of order $O_p(n^{-1/2})$. Finally, it is assumed that the joint cumulants of $nb_a, nb_{ab}, l_r, l_{rs}$, and so forth are of order $O(n)$.

DiCiccio \& Stern (1994b) showed that
\begin{align*}
\bar W(\psi)&=W(\psi)-2\lambda^{ar}\beta_al_r-2\lambda^{ar}b_al_r+2\lambda^{ar}\lambda^{st}\beta_al_{rs}l_t-\lambda^{ar}
\lambda^{su}\lambda^{tv}\beta_a\lambda_{rst} l_u l_v\\
&\hspace{20pt}+\lambda^{ar}\lambda^{bs}\beta_{ab}l_r l_s-\lambda^{ab}\beta_a\beta_b+O_p(n^{-3/2}),
\end{align*}
and it follows that
\begin{align*}
E\{\bar W(\psi)\}&=E\{W(\psi)\}-2\lambda^{ar}E(b_a l_r)+\lambda^{ar}\lambda^{st}\beta_a(2\lambda_{rs,t}+\lambda_{rst})-\lambda^{ab}(\beta_{ab}+\beta_a\beta_b)
+O(n^{-3/2})\\
&=E\{W(\psi)\}+\lambda^{ar}\lambda^{st}\beta_a(2\lambda_{rs,t}+\lambda_{rst})-2\lambda^{ar}\beta_{a/r}
+\lambda^{ab}(\beta_{ab}-\beta_a\beta_b)+O(n^{-3/2}),
\end{align*}
where $\beta_{a/r}=\partial\beta_a/\partial\theta^r$.
For calculating $E\{\bar W(\psi)\}$, we assume that $B(\psi)$ is a function of $Y$ and $\psi$ only, so, in particular, it does not depend on $\phi$. Thus, differentiation of the identity $\beta_a=E\{B_a(\psi)\}$ yields $\beta_{a/b}=E(b_a l_b)+\beta_{ab}$ and $\beta_{a/i}=E(b_a l_i)$ for $i=q+1,\ldots, d$.
It follows that $\lambda^{ar}E(b_a l_r)=\lambda^{ar}\beta_{a/r}-\lambda^{ab}\beta_{ab}$.

To calculate $\mathring E\{W(\psi)\}$, some care is required about the conditional properties of $B_a(\psi)$, $B_{ab}(\psi)$, and so forth.
The quantities $\mathring\beta_a=\mathring E\{B_a(\psi)\}$, $\mathring\beta_{ab}=\mathring E\{B_{ab}(\psi)\}$, etc. are assumed to be of order $O_p(1)$, while
$\mathring b_a=B_a(\psi)-\mathring \beta_a$, $\mathring b_{ab}=B_{ab}(\psi)-\mathring\beta_{ab}$, etc. are assumed to be of order $O_p(n^{-1/2})$.
Finally, it is assumed that the joint conditional cumulants of $n\mathring b_a$, $n\mathring b_{ab}$, $\mathring l_r$, $\mathring l_{rs}$, and so forth are of order $O_p(n)$.

Under the preceding assumptions,
it is possible to determine the orders of the differences $\mathring\beta_a -\beta_a$ and $\mathring \beta_{ab}-\beta_{ab}$. Since $E(\mathring \beta_a)=E[\mathring E\{B_a(\psi)\}]=E\{B_a(\psi)\}=\beta_a$
and ${\rm var}(\mathring\beta_a)={\rm var}[\mathring E\{B_a(\psi)\}]={\rm var}\{B_a(\psi)\}-E[\mathring{\rm var}\{B_a(\psi)\}]=O(n^{-1})-E\{\mathring{\rm var}(\mathring b_a)\}=O(n^{-1})-E\{O_p(n^{-1})\}=O(n^{-1})$,
it follows that $\mathring\beta_a=\beta_a+O_p(n^{-1/2})$. A similar argument shows that $\mathring\beta_{ab}=\beta_{ab}+O_p(n^{-1/2})$.
We assume that differentiation of the identity $\mathring\beta_a=\beta_a+O_p(n^{-1/2})$ yields $\mathring\beta_{a/r}=\beta_{a/r}+O_p(n^{-1/2})$.

Now, define $\mathring\delta_a=\mathring\beta_a-\beta_a$, so that $\mathring\delta_a$ is a function of $\theta$ and $A$ of order $O_p(n^{-1/2})$. Furthermore, $b_a=B_a(\psi)-\beta_a=B_a(\psi)-\mathring\beta_a+\mathring\delta_a=\mathring b_a+\mathring\delta_a$.
To calculate $\mathring E\{\bar W(\psi)\}$, we observe that
\begin{align*}
\bar W(\psi)&=W(\psi)-2\lambda^{ar}\beta_al_r-2\lambda^{ar}b_al_r+2\lambda^{ar}\lambda^{st}\beta_al_{rs}l_t-
\lambda^{ar}\lambda^{su}\lambda^{tv}\beta_a\lambda_{rst} l_u l_v\\
&\hspace{20pt}+\lambda^{ar}\lambda^{bs}\beta_{ab}l_r l_s-\lambda^{ab}\beta_a\beta_b+O_p(n^{-3/2})\\
&=W(\psi)-2\lambda^{ar}\beta_a\mathring l_r-2\lambda^{ar}(\mathring b_a+\mathring\delta_a)\mathring l_r+2\lambda^{ar}\lambda^{st}\beta_a(\mathring l_{rs}+\mathring\Delta_{rs})\mathring l_t-\lambda^{ar}\lambda^{su}\lambda^{tv}\beta_a\lambda_{rst}\mathring l_u\mathring l_v\\
&\hspace{20pt}+\lambda^{ar}\lambda^{bs}\beta_{ab}\mathring l_r \mathring l_s-\lambda^{ab}\beta_a\beta_b+O_p(n^{-3/2}),
\end{align*}
and thus
\begin{align*}
\mathring E\{\bar W(\psi)\}
&=\mathring E\{W(\psi)\}-2\lambda^{ar}\mathring b_a\mathring l_r+2\lambda^{ar}\lambda^{st}\beta_a\mathring\lambda_{rs,t}+\lambda^{ar}\lambda^{su}\lambda^{tv}\beta_a
\lambda_{rst}\mathring\lambda_{uv}\\
&\hspace{20pt}-\lambda^{ar}\lambda^{bs}\beta_{ab}\mathring\lambda_{rs}-\lambda^{ab}\beta_a\beta_b+O_p(n^{-3/2}).
\end{align*}
Barndorff-Nielsen \& Cox (1984) showed that $\mathring E\{W(\psi)\}=E\{W(\psi)\}+O_p(n^{-3/2})$ ; recall that
$\mathring\lambda_{rs}=\lambda_{rs}+O_p(n^{1/2})$ and $\mathring\lambda_{rs,t}=\lambda_{rs,t}+O_p(n^{1/2})$. Then, $\lambda^{ru}\lambda^{st}\mathring\lambda_{ut}=\lambda^{rs}+O_p(n^{-3/2})$, and
\[
\mathring E\{\bar W(\psi)\}
= E\{W(\psi)\}
+\lambda^{ar}\lambda^{st}\beta_a(2\lambda_{rs,t}+\lambda_{rst})-2\lambda^{ar}\mathring E(\mathring b_a \mathring l_r)
-\lambda^{ab}(\beta_{ab}+\beta_a\beta_b)
+O_p(n^{-3/2}).
\]
Now, using the result that $\lambda^{ar}\mathring E(\mathring b_a \mathring l_r)=\lambda^{ar}\mathring\beta_{a/r}-\lambda^{ab}\mathring \beta_{ab}=\lambda^{ar}\beta_{a/r}-\lambda^{ab}\beta_{ab}+O_p(n^{-3/2})$, which holds since $\mathring\beta_{a/r}=\beta_{a/r}+O_p(n^{-1/2})$ and $\mathring\beta_{ab}=\beta_{ab}+O_p(n^{-1/2})$, we have
\begin{align*}
\mathring E\{\bar W(\psi)\}&=E\{W(\psi)\}+\lambda^{ar}\lambda^{st}\beta_a(2\lambda_{rs,t}+\lambda_{rst})-2\lambda^{ar}\beta_{a/r}+
\lambda^{ab}(\beta_{ab}-\beta_a\beta_b)+O_p(n^{-3/2})\\
&=E\{\bar W(\psi)\}+O_p(n^{-3/2}),
\end{align*}
as required. \endex
\par

\section*{References}
\begin{enumerate}
\item
Barndorff-Nielsen, O.\ E.\ (1983).\ On a formula for the conditional distribution of the maximum likelihood estimator.\ {\it  Biometrika} {\bf 70}, 343--65.
\item
Barndorff-Nielsen, O.\ E.\ (1986).\ Inference on full or partial parameters based on the standardized signed log likelihood ratio.\ {\it Biometrika} {\bf 73}, 307--22.
\item
Barndorff-Nielsen, O.\ E.\ and Chamberlin, S.\ R.\  (1994).\ Stable and invariant adjusted directed likelihoods.\ {\it Biometrika} {\bf 81}, 485--99.
\item
Barndorff-Nielsen, O.\ E.\ and Cox, D.\ R.\  (1984).\ Bartlett adjustments to the likelihood ratio statistic and the distribution of the maximum likelihood estimator.\ {\it J.R. Statist. Soc}.\ B {\bf 46}, 483--95.
\item
Barndorff-Nielsen, O.\ E.\ and Cox, D.\ R.\  (1994).\ {\it Inference and
Asymptotics}.\ Chapman \& Hall, London.
\item
Barndorff-Nielsen, O.\ E.\ and Hall, P.\   (1988).\ On the level-error
after Bartlett adjustment of the likelihood ratio statistic.\ {\it Biometrika}
{\bf 75}, 374--8.
\item
Bickel, J.\ K.\ and Ghosh, J.\ K.\  (1990).\ A decomposition for the likelihood ratio statistic and the Bartlett correction -- a Bayesian argument.\ {\it Ann.\ Statist.\ } {\bf 18}, 1070--90.
\item
Cordeiro, G. and Ferrari, S.\ L.\ de P.\  (1991).\ A modified score test statistic having chi-squared distribution to order $n^{-1}$.\ {\it Biometrika} {\bf 78}, 573--82.
\item
Cox, D.\ R.\  (1980).\
Local ancillarity.\
{\it Biometrika} {\bf 67}, 279--86.
\item
Cox, D.\ R.\ and Reid, N.\  (1987).\  Parameter orthogonality and approximate conditional inference
(with discussion).\ {\it J.R. Statist. Soc}.\ B {\bf 53}, 79--109.
\item
DiCiccio, T.\ J., Kuffner, T.\ A. and Young, G. \ A. \ (2014). Inferential correctness and the parametric bootstrap. \ In preparation.
\item
DiCiccio, T.\ J.\ and Martin, M.\ A.\  (1993).\ Simple modifications for signed roots of likelihood
ratio statistics.\  {\it J.R. Statist. Soc}.\ B {\bf 55}, 305--16.
\item
DiCiccio, T.\ J.\ and Stern, S.\ E.\  (1994a).\ Constructing approximately standard normal pivots from signed roots of adjusted likelihood ratio statistics.\
{\it Scand. J. Statist.}, {\bf 21}, 447--60.
\item
DiCiccio, T.\ J.\ and Stern, S.\ E.\  (1994b).\ Frequentist and Bayesian Bartlett correction of test statistics based on adjusted profile likelihoods.\ {\it J.R. Statist. Soc}.\ B {\bf 56}, 397--408.
\item
DiCiccio, T.\ J.\ and Young, G.\ A.\  (2008).\ Conditional properties
of unconditional parametric bootstrap procedures for inference in exponential
families.\ {\it Biometrika} {\bf 95}, 747--58.
\item
DiCiccio, T.\ J.\ and Young, G.\ A.\  (2010).\
Computer-intensive conditional inference. In {\it Complex Data Modeling and Computationally Intensive Statistical Methods} (P. Mantovan, P. Secchi eds.), 137--150. Springer-Verlag Italia, Milan.
\item
DiCiccio, T.\ J., Martin, M.\ A.\ and Stern, S.\ E.\  (2001).\  Simple and accurate one-sided inference from signed roots of likelihood ratios.\ {\it Can.\ J.\ Statist}.\ {\bf 29}, 67--76.
\item
DiCiccio, T.\ J., Martin, M.\ A., Stern, S.\ E., and Young, G.\ A.\  (1996).\
Information bias and adjusted
profile likelihood.\ {\it J.R. Statist. Soc}.\ B {\bf 58}, 189--203.
\item
Efron, B.\  (1975).\ Defining the curvature of a statistical problem (with applications to second order efficiency) (with discussion).\ {\it Ann.\ Statist.\ } {\bf 3}, 1189--1242.
\item
Efron, B.\ and Hinkley, D.\ V.\  (1978).\ Assessing the accuracy of the maximum likelihood estimator: Observed versus expected Fisher information (with discussion).\ {\it Biometrika} {\bf 65}, 457--87.
\item
Firth, D.\  (1993).\ Bias reduction of maximum likelihood estimates.\ {\it Biometrika} {\bf 80}, 27--38.
\item
Harris, P.\  (1985).\ An asymptotic expansion for the null distribution of the efficient score statistic.\ {\it Biometrika} {\bf 72}, 653--9.
\item
Lawley, D.\ N.\  (1956).\ A general method for approximating to the distribution of likelihood ratio criteria.\ {\it  Biometrika} {\bf 43}, 295--303.
\item
Lee, S.\ M.\ S.\ and Young, G.\ A.\  (2005).\ Parametric bootstrapping
with nuisance parameters.\ {\it Stat. Prob. Letters} {\bf 71}, 143--53.
\item
McCullagh, P.\  (1984).\ Local sufficiency.\ {\it Biometrika} {\bf 71}, 233--44.
\item
McCullagh, P.\  (1987).\ {\it Tensor Methods in Statistics}.\ Chapman \& Hall, London.
\item
McCullagh, P.\  (1992).\ Conditional inference and Cauchy models.\  {\it Biometrika} {\bf 79}, 247--59.
\item
McCullagh, P.\ and Tibshirani, R.\  (1990).\ A simple method for the adjustment of profile likelihoods.\ {\it J. Roy. Statist. Soc. B} {\bf 52}, 325--44.
\item
Pace, L.\  and Salvan, A.\  (1994).\ The geometric structure of the expected/observed likelihood expansions.\
{\it Ann. Inst. Stat. Math.} {\bf 46}, 649--66.
\item
Pace, L.\ and Salvan, A.\  (1997).\ {\it Principles of Statistical Inference:
from a Neo-Fisherian Perspective}.\ World Scientific, Singapore.
\item
Pierce, D.\ A.\  (1975).\ Discussion of paper by B. Efron.\ {\it Ann.\ Statist.\ } {\bf 3}, 1219--21.
\item
Pierce, D.\ A.\ and Bellio, R.\  (2006).\ Effects of the reference set on frequentist inferences.\
{\it Biometrika} {\bf 93}, 425--38.
\item
Severini, T.\ A.\  (1990).\ Conditional properties of likelihood-based significance tests.\ {\it Biometrika} {\bf 77}, 343--52.
\item
Severini, T.\ A.\  (1998).\ An approximation to the modified profile likelihood function.\ {\it Biometrika} {\bf 85}, 403--11.
\item
Severini, T.\ A.\  (2000).\ {\it Likelihood methods in Statistics}.\ Oxford University Press, Oxford.
\item 
Skovgaard, I.\ M.\  (1996).\ An explicit large-deviation approximation to one-parameter tests.\
{\it Bernoulli} {\bf 2}, 145-65.
\item
Tierney, L.\ and Kadane, J.\ B.\  (1986).\
Accurate approximations for posterior moments and marginal densities.\
{\it J. Amer. Statist. Assoc.} {\bf 81}, 82--6.
\item Young, G.\ A.\  (2009).\ Routes to higher-order accuracy in parametric inference.\ {\it Aust. N. Z. J. Stat.\ } {\bf 51}, 115--26.
\item Young, G.\ A.\ and Smith, R.\ L.\  (2005).\ {\it Essentials of Statistical
Inference}.\ Cambridge University Press, Cambridge.
\end{enumerate}

\end{document}